\documentclass[nonblindrev]{informs3}

\OneAndAHalfSpacedXI 

\usepackage{endnotes}
\let\footnote=\endnote

\usepackage{natbib}
 \bibpunct[, ]{(}{)}{,}{a}{}{,}%
 \def\bibfont{\small}%
\TheoremsNumberedThrough     
\ECRepeatTheorems

\EquationsNumberedThrough    

\usepackage{hyperref}
\hypersetup{
     colorlinks = true,
     linkcolor = blue,
     anchorcolor = blue,
     citecolor = blue,
     filecolor = blue,
     urlcolor = blue
     }
\usepackage{color}
\usepackage{amsmath}
\usepackage{amssymb}
\usepackage{amsfonts}
\usepackage[bb=boondox]{mathalfa}
\usepackage{mathtools}
\usepackage{epsfig}
\usepackage{graphicx}
\usepackage{graphics}
\usepackage{xcolor}
\usepackage{float}
\usepackage[shortlabels]{enumitem}
\usepackage{multirow}
\usepackage{caption}
\usepackage{subcaption}
\usepackage{xspace}
\usepackage{booktabs}

\usepackage{tikz}
\usepackage{pgfplots}
\pgfplotsset{compat=1.17}
\usetikzlibrary{arrows.meta}
\usepackage{marvosym}
\usetikzlibrary{positioning,fit,arrows.meta,backgrounds}
\usetikzlibrary{patterns,shapes,positioning}
\usetikzlibrary{decorations.markings}
\usepgfplotslibrary{fillbetween}
\usetikzlibrary{decorations.text}
\usetikzlibrary{arrows, decorations.pathmorphing}
\usepgfplotslibrary{fillbetween}
\usetikzlibrary{patterns,shapes,positioning}

\newcommand{\minimize}[1]{\underset{#1}{\text{minimize}}}
\newcommand{\maximize}[1]{\underset{#1}{\text{maximize}}}
\newcommand{\st}{\text{subject to}}
\newcommand{\bs}[1]{\boldsymbol{#1}}

\newcommand\norm[1]{\left\lVert#1\right\rVert}
\newcommand*{\vertbar}{\textcolor{gray}{\rule[0.4ex]{0.5pt}{0.2ex}}}

\definecolor{maincolor}{RGB}{0,0,0}
\definecolor{red}{HTML}{e05a87} 

\usepackage{algorithm2e}
\usepackage{xcolor}

\newcommand{\xvbox}[2]{\makebox[#1][l]{#2}}

\SetAlFnt{\footnotesize}

\SetAlCapSty{xAlCapSty}

\SetCommentSty{xCommentSty}

\SetNlSty{mynlfont}{}{} 

\LinesNumbered

\SetSideCommentRight

\DontPrintSemicolon

\RestyleAlgo{algoruled}

\begin{document}


\RUNAUTHOR{Dvorkin et al.}

\RUNTITLE{Privacy-Preserving Convex Optimization}

\TITLE{Privacy-Preserving Convex Optimization: When Differential Privacy Meets Stochastic Programming}

\ARTICLEAUTHORS{%
\AUTHOR{Vladimir Dvorkin}
\AFF{Department of Electrical Engineering and Computer Science, University of Michigan, Ann Arbor, MI 48103, \EMAIL{dvorkin@umich.edu}} 
\AUTHOR{Ferdinando Fioretto}
\AFF{Department of Computer Science, University of Virginia, Charlottesville, VA 22903, \EMAIL{ fioretto@virginia.edu}}
\AUTHOR{Pascal Van Hentenryck}
\AFF{H. Milton Stewart School of Industrial and Systems Engineering, Georgia Institute of Technology, Atlanta, GA 30332, \EMAIL{pascal.vanhentenryck@isye.gatech.edu}}
\AUTHOR{Pierre Pinson}
\AFF{Dyson School of Design Engineering, Imperial College London, London, SW7 2AZ, U.K, \EMAIL{p.pinson@imperial.ac.uk}}
\AUTHOR{Jalal Kazempour}
\AFF{Department of Wind and Energy Systems, Technical University of Denmark, Lyngby, 2800, Denmark, \EMAIL{jalal@dtu.dk}}
} 

\ABSTRACT{%
Convex optimization finds many real-life applications, where -- optimized on real data -- optimization results may expose private data attributes (e.g., individual health records, commercial information), thus leading to privacy breaches. To avoid these breaches and formally guarantee privacy to optimization data owners, we develop a new privacy-preserving perturbation strategy for convex optimization programs by combining stochastic (chance-constrained) programming and differential privacy. Unlike standard noise-additive strategies, which perturb either optimization data or optimization results, we express the optimization variables as functions of the random perturbation using linear decision rules; we then optimize these rules to accommodate the perturbation within the problem's feasible region by enforcing chance constraints. This way, the perturbation is feasible and makes different, yet adjacent in the sense of a given distance function, optimization datasets statistically similar in randomized optimization results, thereby enabling probabilistic differential privacy guarantees. The chance-constrained optimization additionally internalizes the conditional value-at-risk measure to model the tolerance towards the worst-case realizations of the optimality loss w.r.t. the non-private solution. We demonstrate the privacy properties of our perturbation strategy analytically and through optimization and machine learning applications. 
}%




\maketitle

%

\vspace{-0.5cm}
\section{Introduction}

As data become easier to extract and improve in quality, designing optimization and machine learning models that can leverage these data has become a focal point in finance, healthcare, energy, transportation, and other fields. Often data include private and sensitive records, such as personalized health and financial records, which could be revealed when releasing model parameters or results of model applications, thus leading to privacy breaches for data owners.  

The privacy breaches in linear programming originate from inverse optimization techniques that reveal unknown parameters of an optimization program from its optimal solution \citep{ahuja2001inverse}. 
Their applications to electricity markets reveal classified transmission parameters \citep{birge2017inverse} and private bids of market participants \citep{ruiz2013revealing,mitridati2017bayesian}. Such information weaponizes strategic market actors that maximize the profits at the expense of deteriorating social welfare \citep{chen2019learning}. Privacy breaches are relevant for machine learning models represented by convex, quadratic programs. While model parameters are seen as some data aggregation, they are still a function of the underlying training datasets. In prediction problems, linear regression weights are unique on a particular dataset. In classification problems, the parameters of support vector machines (SVM), one of the most robust classifiers \citep{xu2009robustness}, are unique w.r.t. the marginal data points \citep{burges1999uniqueness}. These properties are exploited by reconstruction attacks that outcome one or more training samples and their respective training labels, or by membership attacks that determine the presence of a particular sample in a training dataset \citep{rigaki2020survey}. The privacy risks are also seen in semidefinite programming, e.g., for maximum-volume inscribed ellipsoid computation for reducing the conservatism of a robust uncertainty set, employed -- to name one example -- in online pricing algorithms \citep{cohen2020feature}. Ellipsoid parameters together with all but one element of an uncertainty set reveal the remaining one element, hence violating privacy. 

Providing datasets with formal privacy guarantees  is accomplished with controlled perturbations of computations to decouple datasets from computational outcomes. Differential privacy provides such guarantees by calibrating perturbations to the sensitivity of computational outcomes to changes in datasets , so as to obfuscate these changes \citep{dwork2006calibrating}. This makes different, yet adjacent in the sense of a prescribed distance function, datasets  statistically similar in randomized outcomes, thereby providing a probabilistic privacy guarantee. At first, differentially private computations have guaranteed individual privacy, i.e., masking the presence of an individual record or its attributes, and later a metric-based differential privacy has been introduced by \citet{chatzikokolakis2013broadening} to obfuscate arbitrary distances between datasets, enabling privacy of data subsets or entire datasets, which is relevant for many applications of convex optimization. In Section \ref{subsec:related_work} we briefly review differentially private optimization, and in Section \ref{subsec:contribution} we state the contribution w.r.t. the state-of-the-art methods. 

\subsection{Differentially Private Convex Optimization}\label{subsec:related_work}

Differential privacy for optimization datasets has been provided using different perturbation strategies. \citet{chaudhuri2011differentially} developed private learning algorithms that perturb either the parameters of optimization (objective perturbation, also referred to as {\it input} perturbation) or the results of optimization ({\it output} perturbation), providing differentially private logistic regression and SVM models. Privacy leakages in learning problems are also limited by iterative randomized algorithms, such as noisy gradient descent algorithms by \citet{bassily2014private}, that enjoy gradient perturbation at every iteration. The advantage of \citep{bassily2014private} over \citep{chaudhuri2011differentially} is that the sensitivity is easier to bound by hyper-parameters, i.e., by a gradient clipping constant. However, unlike a single-shot optimization in \citep{chaudhuri2011differentially}, the privacy budget in \citep{bassily2014private} accumulates over iterations. To void this issue, composition theorems of differential privacy apply to scale perturbation parameters w.r.t. the maximum number of iterations \citep{dwork2014algorithmic}, yet leading to a large variance of the result. Building on ideas of input, output and gradient perturbations, differentially private variants of regression and classification models also appear in \citep{talwar2015nearly,wang2017differentially,wang2018empirical,sheffet2019old,sun2021dpwss}.

Guaranteeing privacy for constrained programs is relevant for optimization problems whose constraints contain private information, such as resource allocation problems. Unlike learning problems, they are not restricted to strictly monotone objective functions and the standing assumptions of \citet{chaudhuri2011differentially} are not necessarily met. The privacy in linear allocation problems has been thus addressed using algorithmic solutions, such as ones in \citep{hsu2014privately} based on privacy-preserving variants of the dense multiplicative weights algorithm. The algorithm, however, only returns an approximate solution. \textcolor{maincolor}{In many applications of interest, such as constrained resource allocation, the approximate solution may not satisfy the constraints, prompting the decision-maker to possibly compromise on privacy in favor of feasibility.} This problem has been partially addressed by \citet{munoz2021private} with a truncated dataset perturbation which ensures that the perturbed solution remains in the feasible region. Truncated  perturbations, however, do not apply to optimization problems that include equality constraints. 

An alternative solution to the privacy of resource allocation data is to decompose a global optimization problem into a set of sub-problems that restore the solution to the global problem over iterations with limited information exchange, e.g., using the alternating direction method of multipliers, ADMM \citep{boyd2011distributed}. Although datasets are distributed among sub-problems and not explicitly shared, they tend to leak through the exchanged ADMM variables \citep{dvorkin2020differentially}. To void this leakage, differentially private ADMM variants -- that perturbed either primal or dual variables -- have been introduced in \citep{han2016differentially,dvorkin2020differentially,cao2020differentially}. Even though ADMM convergence is robust to minor perturbations \citep{majzoobi2019analysis}, it may fail to converge when privacy guarantees extend beyond one iteration \citep{dvorkin2020differentially}. Moreover, although the design of private distributed algorithms ensures the per-iteration privacy of distributed datasets, it does not guarantee the privacy of the global mapping from the entire dataset to optimization results. Similar limitations of the private ADMM are also shown in machine learning applications \citep{zhang2018recycled,huang2019dp}.

\subsection{Contributions}\label{subsec:contribution}
The summary in Section \ref{subsec:related_work} highlights that there is no general framework to provide differential privacy to arbitrary convex optimization programs, e.g., the algorithms for differentially private SVM due to \citep{chaudhuri2011differentially} do not apply to private linear programming in \citep{munoz2021private}, and vice versa. The main contribution of this work is to develop a general framework for convex optimization programs with a formal privacy guarantee. Without making assumptions on problem type, the framework addresses a class of {\it conic optimization} that encompasses linear, quadratic, semidefinite and other convex program types. For a deterministic conic program, we provide a stochastic optimization counterpart, which maps optimization data to solution space {\it randomly} in a way that adjacent datasets produce statistically similar results up to prescribed differential privacy requirements. However, unlike optimization-agnostic input and output perturbation strategies, this strategy respects {\it feasibility} and {\it optimality} criteria internalized into the stochastic formulation. More formally, we make the following technical contributions:
\begin{enumerate}[topsep=0pt]\itemsep0em 
    \item We provide differential privacy to datasets of convex optimization programs by enforcing a linear dependency of optimization variables on the privacy-preserving perturbation. This dependency is enforced using linear decision rules, popularized in the decision-making literature to accommodate the uncertainty of datasets \citep{ben2004adjustable,shapiro2005complexity,kuhn2011primal}. In the present case, however, they are used to randomize the solution of a deterministic program, so as to provide privacy to the program's dataset. The optimized linear decision rules enable private linear queries, such as identity query to release the entire solution vector (e.g., parameters of machine learning models), or sum query to release the optimal objective function value (e.g., allocation cost).
    \item We leverage a chance-constrained linear decision rule optimization to impose feasibility requirements on differentially private optimization results. Similar to decision-making under uncertainty, we identify a trade-off between feasibility and the expected cost of the privacy-preserving solution. We also establish a monotonic relationship between differential privacy parameters and the cost of privacy in our stochastic framework.  
    \item The proposed framework evaluates the expected suboptimality of the privacy-preserving solution as a distance between the non-private, deterministic solution and the chance-constrained solution. For linear queries made on strict subsets of the solution vector, the chance-constrained optimization extends to optimize the worst-case optimality loss w.r.t. a chosen risk measure, e.g., conditional value--at--risk \citep{rockafellar2000optimization}, enabling additional trade-offs between the expected and worst-case solution suboptimality due to perturbation. 
\end{enumerate}

We call our framework {\it program perturbation}, as it requires solving a stochastic variant of the original program, although preserving its data and structure. For unconstrained optimization (e.g., most of learning problems), input and output perturbations are competitive with the program perturbation when the parameters of perturbation are marginal, e.g., when protecting a single entry of a training dataset. The benefit of program perturbation reveals when perturbation magnitude is significantly far from marginal, e.g., when the distance between datasets, that must be made indistinguishable, is large. This is the case, for example, of a regression model that explains a physical phenomenon on heterogeneous datasets, but the perturbation of the regression model w.r.t. the distance between these datasets leads to physically inconsistent results. For constrained optimization (e.g., linear programming), program perturbation systematically outperforms either input or output perturbation regardless of perturbation parameters. 

\subsection{Paper Organization}

The remainder of this paper is organized as follows. Section \ref{sec:preliminaries} provides preliminaries on conic optimization and differential privacy, as well as states the problem using an example. Section \ref{sec:main} explains how perturbations are internalized using linear decision rules and chance-constrained optimization in risk-neutral and risk-averse cases. Section \ref{sec:application} demonstrates applications to a range of problems from the electric power system domain. Section \ref{sec:conclusion} concludes, while proofs and additional modeling details are relegated to the e-companion. We also release all data, codes and tutorials to replicate our results in online repository in \citep{dvorkin2022github}. 

\paragraph{Notation:} lower- and upper-case letters denote column vectors and matrices, respectively. $\mathbb{0}$ and $\mathbb{1}$ respectively denote all-zero and all-one vectors (or matrices), and $I$ is an identity matrix; they are appended with a subscript whenever dimension is ambiguous.  Operation $\circ$ is the element-wise product. Operator $\text{diag}[x]$ returns an $n\times n$ diagonal matrix with elements of vector $x\in\mathbb{R}^{n}$, $A_{i}$ returns the $i$\textsuperscript{th} row of matrix $A$ as a vector, and operator $\text{Tr}[A]$ returns the trace of $A$. Symbol $^\top$ is transposition. $\norm{\cdot}_{p}$ denotes a $p-$norm; without subscript $p$, an $\ell_2-$norm is implied.

\section{Preliminaries and Problem Statement}\label{sec:preliminaries}
\subsection{Convex Program as a Map}
We consider the problem of a data curator who owns or collects optimization datasets from other owners and maps them to the solution of a convex optimization program, i.e., by solving
\begin{subequations}\label{problem:base}
\begin{align}
    x(\mathcal{D}) = \underset{x}{\text{arg\!\;min}}\quad& c^{\top}x\label{problem_base_obj}\\
    \st\quad&
    b - Ax \in\mathcal{K},\label{problem_base_con}
\end{align}
\end{subequations}
where $x\in\mathbb{R}^{n}$ is the vector of optimization variables,  parameters $A\in\mathbb{R}^{m\times n}$, $b\in\mathbb{R}^{m}$ and $c\in\mathbb{R}^{n}$ are optimization data, and $\mathcal{K}\subseteq \mathbb{R}^{m}$ denotes an $m-$dimensional cone. Notation $x(\mathcal{D})$ denotes the map from datasets to the optimal solution $x^{\star}\in\mathbb{R}^{n}$ of program \eqref{problem:base} that minimizes the objective function given all feasible choices in \eqref{problem_base_con}. By optimization dataset, we understand a tuple $\mathcal{D}=( A,b,c)$ with $n+m(1+n)$ real values. Let all optimization datasets belong to a finite domain of real values $\mathbb{D}$, called dataset \textit{universe}. In the following, we are interested in datasets $\mathcal{D}\in\mathbb{D}$ that ensure solution existence. \textcolor{maincolor}{We are also interested in optimization maps $x(\mathcal{D})$ that yield a unique solution on a particular dataset $\mathcal{D}$, as they present the most interesting cases as shown below in Figure \ref{fig:examples}.}

\textcolor{maincolor}{This conic formulation is very general and many convex problems, including linear, quadratic and semidefinite programs, can be represented by \eqref{problem:base} and solved by primal-dual interior-point methods \citep{nesterov1998primal,dahl2022primal}.} For example, to solve a linear program, the cone $\mathcal{K}=\mathbb{R}_{+}^{n}$ must be a non-negative orthant. A linear regression problem that minimizes the norm $\norm{Fy-g}$ by choosing regression weights $y\in\mathbb{R}^{k}$ is also representable by formulation \eqref{problem:base} when dataset $\mathcal{D}=\left(\begin{bsmallmatrix}
-1 & \\& -F
\end{bsmallmatrix},
\begin{bsmallmatrix}
\textcolor{white}{-}0 \\-g
\end{bsmallmatrix},
\begin{bsmallmatrix}
1 \\\mathbb{0}
\end{bsmallmatrix}\right),$ 
solution vector $x=
\begin{bsmallmatrix}
t \\y
\end{bsmallmatrix}
$  and cone
$\mathcal{K}=\mathcal{Q}^{k+1}$ 
is a second-order cone, where auxiliary variable $t$ computes the fitting loss. 
We refer the reader to e-companion \ref{app:conic_representations} for the representation of problems of interest by program \eqref{problem:base}. 

In our setting, the curator answers linear queries on the solution of program \eqref{problem:base}, i.e., it releases any linear transformation of vector $x^{\star}$. For example, in the linear programming case, the weighted sum query $c^{\top}x^{\star}$ releases the optimal objective function value. The identity query in the linear regression case releases regression weights $y$ and fitting loss $t$ induced on a particular dataset.

\subsection{Differential Privacy}\label{subsec:DP}
\begin{figure}
    \FIGURE{
    \begin{subfigure}[t]{0.28\textwidth}
        \centering
        \includegraphics[width=1\textwidth]{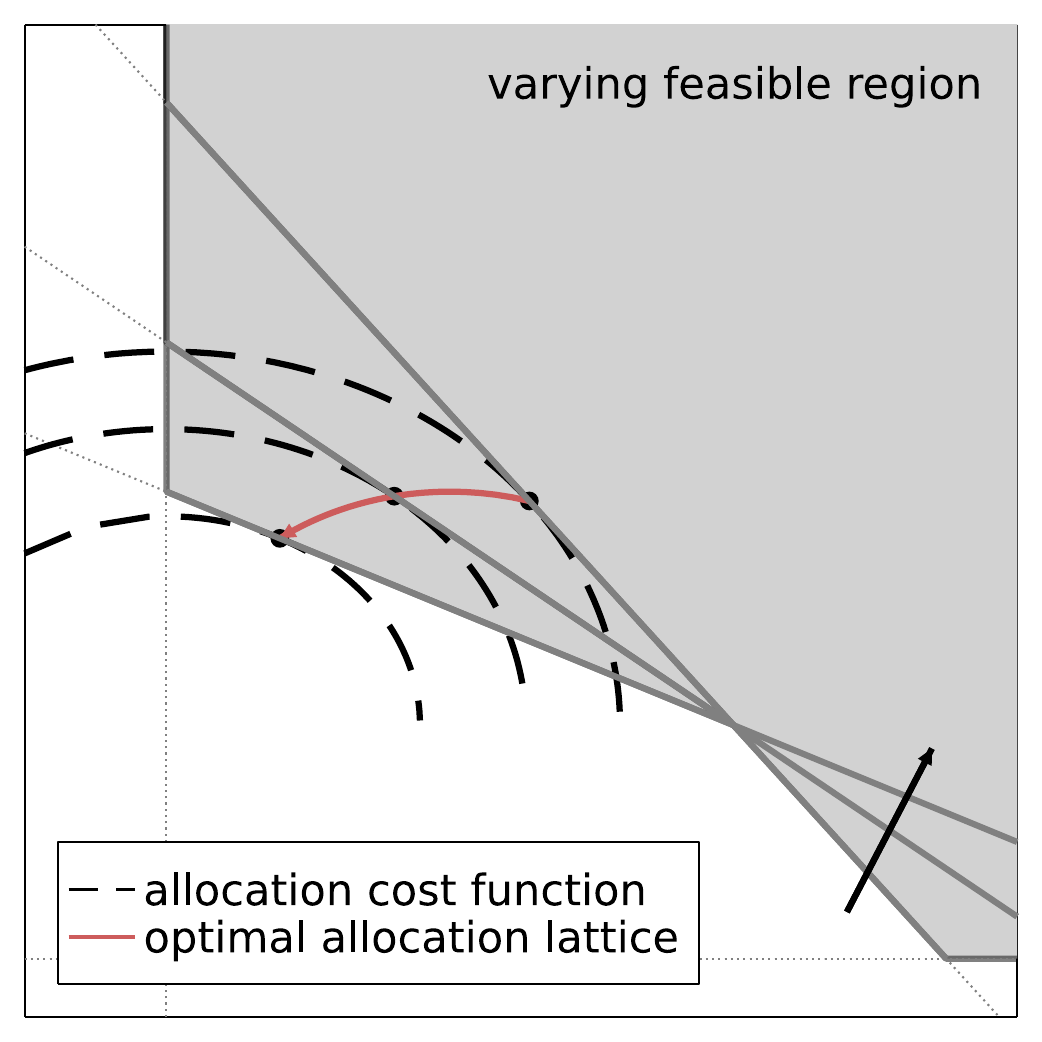}
        \caption{Resource allocation}
    \end{subfigure}%
    ~ \hspace{0.25cm}
    \begin{subfigure}[t]{0.28\textwidth}
        \centering
        \includegraphics[width=1\textwidth]{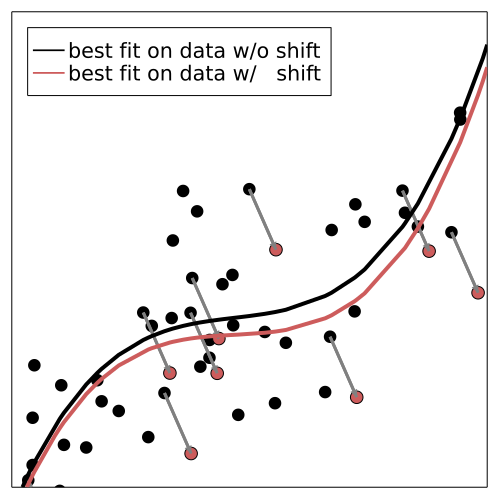}
        \caption{Regression analysis}
    \end{subfigure}%
    ~ \hspace{0.25cm}
    \begin{subfigure}[t]{0.28\textwidth}
        \centering
        \includegraphics[width=1\textwidth]{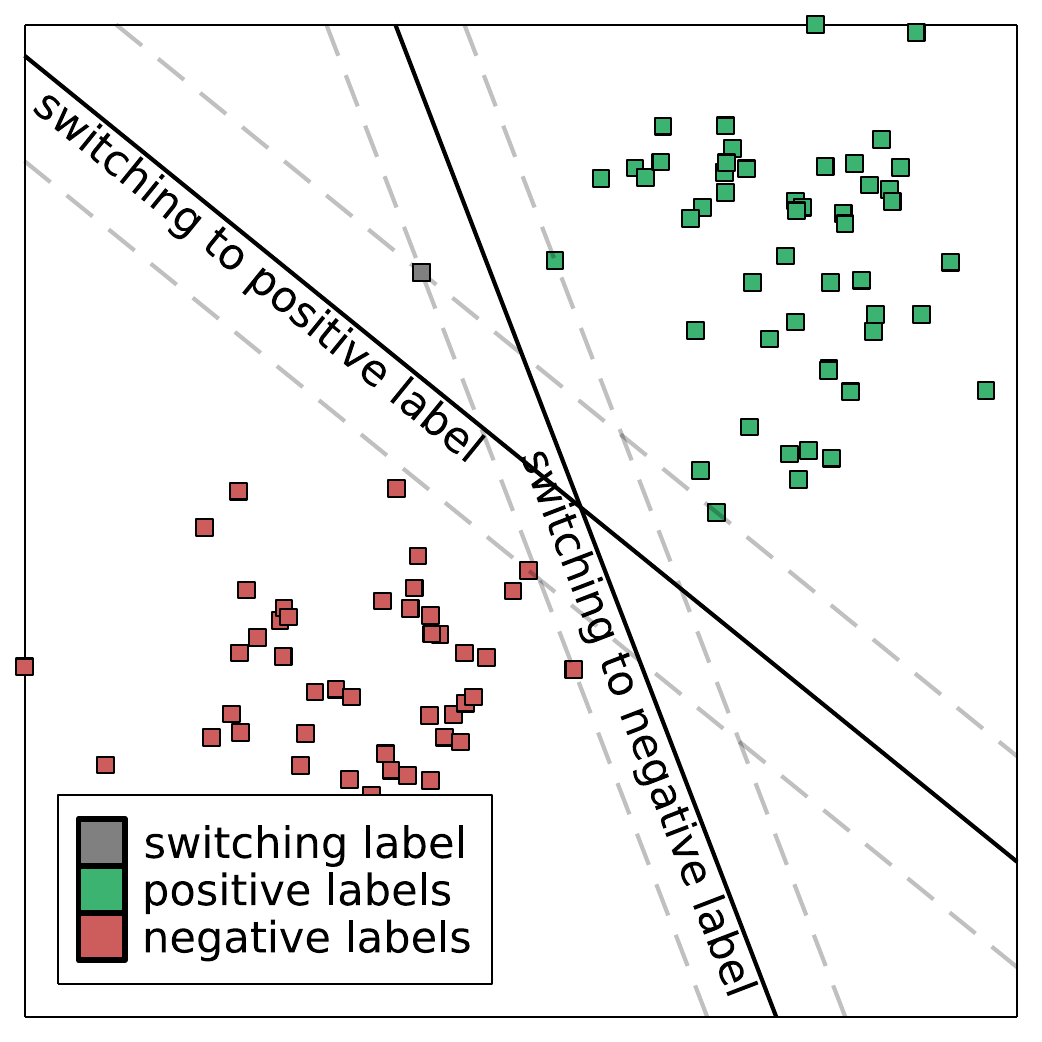}
        \caption{SVM classification}
    \end{subfigure}%
    }
    {Privacy leakages in convex optimization: (a) The optimal allocation always lies at the boundary of the feasible region; hence, changes in the feasible region are directly exposed in the optimal allocation cost. (b) In regression analysis, any data shift is directly exposed in the parameters of the optimal fitting model; here, the shift is applied to a subset of data points, but even a change of a single datum would be exposed. (c) In support vector machines classification, a separating hyperplane is sensitive to marginal data points; here, switching the label of the marginal data point is exposed in hyperplane parameters. We refer the reader to online repository \citep{dvorkin2022github} for more examples. \label{fig:examples}}{}
\end{figure}

When map $x(\mathcal{D})$ yields a unique solution on a dataset $\mathcal{D}$, optimization results on two different datasets \textcolor{maincolor}{may} disclose their underlying differences: one distinguishes between datasets $\mathcal{D},\mathcal{D}'\in\mathbb{D}$ by querying $x(\mathcal{D})$ and $x(\mathcal{D}')$ or their certain linear transformations, as illustrated for resource allocation and machine learning problems in Figure \ref{fig:examples}. \textcolor{maincolor}{This is particular relevant for constraint optimization problems where constraints are built of private data, as those problems studied later in Section \ref{sec:application}.}  Some dataset entries may constitute private information which may unintentionally leak when answering optimization queries.

The differential privacy goal is thus to make adjacent optimization datasets statistically similar when answering optimization queries. By adjacent, we understand any two datasets $\mathcal{D},\mathcal{D}'\in\mathbb{D}$ that are close one to another in the sense of Euclidean distance $\norm{\mathcal{D}-\mathcal{D}'}\leqslant\alpha$ bounded by some prescribed parameter $\alpha\geqslant0;$ \textcolor{maincolor}{sometimes, we will denote adjacency of two datasets by $\sim_{\alpha}$, i.e., $\mathcal{D}\sim_{\alpha}\mathcal{D}'$. This choice is inspired by the metric-based privacy due to \citet{chatzikokolakis2013broadening} and allows to work with different participation serious of individuals in a dataset. For example, we cover cases when adjacent vectors are different only in one entry up to parameter $\alpha$, or may different in several entries as long as the distance between the two vectors is bounded by $\alpha$. The motivation behind this bound is to address an arbitrary large universe $\mathbb{D}$ where there exist datasets that are so far from each other, that hiding the distance between them is not practical.}

Differential privacy guarantees for optimization datasets are achieved through randomization. Let $\tilde{x}(\mathcal{D})$ be a randomized counterpart  of \textcolor{maincolor}{deterministic} optimization map $x(\mathcal{D})$ \textcolor{maincolor}{in the sense} that $\tilde{x}(\mathcal{D})$ maps dataset $\mathcal{D}$ to optimization results randomly. Map $\tilde{x}(\mathcal{D})$ provides differential privacy for $\alpha-$adjacent datasets if it satisfies the following basic definition \citep{dwork2014algorithmic}.  e
\begin{definition}[Differential privacy]\label{def:DP} A random map  $\tilde{x}:\mathbb{D}\mapsto\textcolor{maincolor}{\text{Range}(\tilde{x})}$ is $\varepsilon-$differentially private if  $\forall\hat{x}\in\textcolor{maincolor}{\text{Range}(\tilde{x})}$ and $\forall\mathcal{D},\mathcal{D}'\in\mathbb{D}$ satisfying $\norm{\mathcal{D}-\mathcal{D}'}\leqslant\alpha$ for some $\alpha\geqslant0$, it holds that
\begin{align*}
    \text{Pr}[\tilde{x}(\mathcal{D})\;\textcolor{maincolor}{\in}\;\hat{x}]\leqslant\text{Pr}[\tilde{x}(\mathcal{D}') \;\textcolor{maincolor}{\in}\;\hat{x}]\text{exp}(\varepsilon),
\end{align*}
where probability $\text{Pr}$ is taken over randomness of  $\tilde{x}$.
\end{definition}

According to this definition, the probabilities of observing the same optimization result on adjacent datasets are similar up to prescribed privacy parameter $\varepsilon>0$, \textcolor{maincolor}{termed \textit{privacy budget}}. Consequently, the results of random maps $\tilde{x}(\mathcal{D})$ and $\tilde{x}(\mathcal{D}')$ are more statistically similar when smaller privacy budget. An $\varepsilon-$DP optimization map enjoys the following important privacy property. 
\begin{theorem}[Post-processing immunity]\label{th:post_processing}\normalfont
Let $\tilde{x}$ be an $\varepsilon-$DP map and $g$ be some other, independent map from the set of outcomes of $\tilde{x}$ to an arbitrary set. Then, $g\circ\tilde{x}$ is $\varepsilon-$DP.
\end{theorem}

An optimization map can provide DP to input datasets by perturbing optimization result with a data-independent random noise. To ensure data-independence, the noise is calibrated to the worst-case sensitivity of the map to adjacent datasets \citep{dwork2006calibrating}. In our notation, the worst-case sensitivity is $\Delta_{p}=\text{max}_{\mathcal{D},\mathcal{D}'}\;\norm{x(\mathcal{D}) - x(\mathcal{D}')}_{p},$ for all $(\mathcal{D},\mathcal{D}')\in\mathbb{D}$ satisfying $\norm{\mathcal{D}-\mathcal{D}'}\leqslant\alpha$. 


The {\it Laplace mechanism} acting on $\Delta_{1}-$sensitivity is commonly used to achieve DP. \textcolor{maincolor}{Let $\text{Lap}(s)^{n}$ denote the $n-$dimensional Laplace distribution with zero mean and scale $s$. The Laplace mechanism is defined as follows \citep{dwork2006calibrating}.}
\begin{theorem}[Laplace mechanism]\label{th:Laplace}\normalfont Let $x$ be a map from datasets to $\mathbb{R}^{n}$. The Laplace mechanism $x(\mathcal{D})+\hat{\zeta}$, where perturbation $\hat{\zeta}$ is drawn from $\text{Lap}(\Delta_{1}/\varepsilon)^{n}$, is $\varepsilon-$DP. 
\end{theorem}


We can now formulate two preliminary privacy-preserving strategies for convex optimization. 
\begin{definition}[Output perturbation]\label{def:OP}\normalfont The output perturbation is $x(\mathcal{D})+\hat{\zeta}$ with perturbation $\hat{\zeta}$ from $\text{Lap}(\Delta_{1}/\varepsilon)^{n}$, where $\Delta_{1}$ is the worst-case $\ell_{1}-$sensitivity of the map to adjacent datasets.
\end{definition}
\begin{definition}[Input perturbation]\label{def:IP}\normalfont For dataset universe $\mathbb{D}\subset\mathbb{R}^{k}$, the input perturbation is twofold: data perturbation $\tilde{\mathcal{D}}=\mathcal{D}+\hat{\zeta}$ with perturbation $\hat{\zeta}$ from $\text{Lap}(\alpha/\varepsilon)^{k}$, where $\alpha$ is the worst-case $\ell_{1}-$sensitivity $\Delta_{1}$ of an identity query to $\alpha-$adjacent datasets, is followed by the map $x(\tilde{\mathcal{D}})$.
\end{definition}
\begin{proposition}\label{prop:OP_IP_dp}\normalfont
The output and input perturbation strategies are $\varepsilon-$ differentially private.
\end{proposition}

\subsection{Limits of the Standard Perturbation Strategies}
\begin{figure}
\FIGURE{
\includegraphics[width=0.75\textwidth]{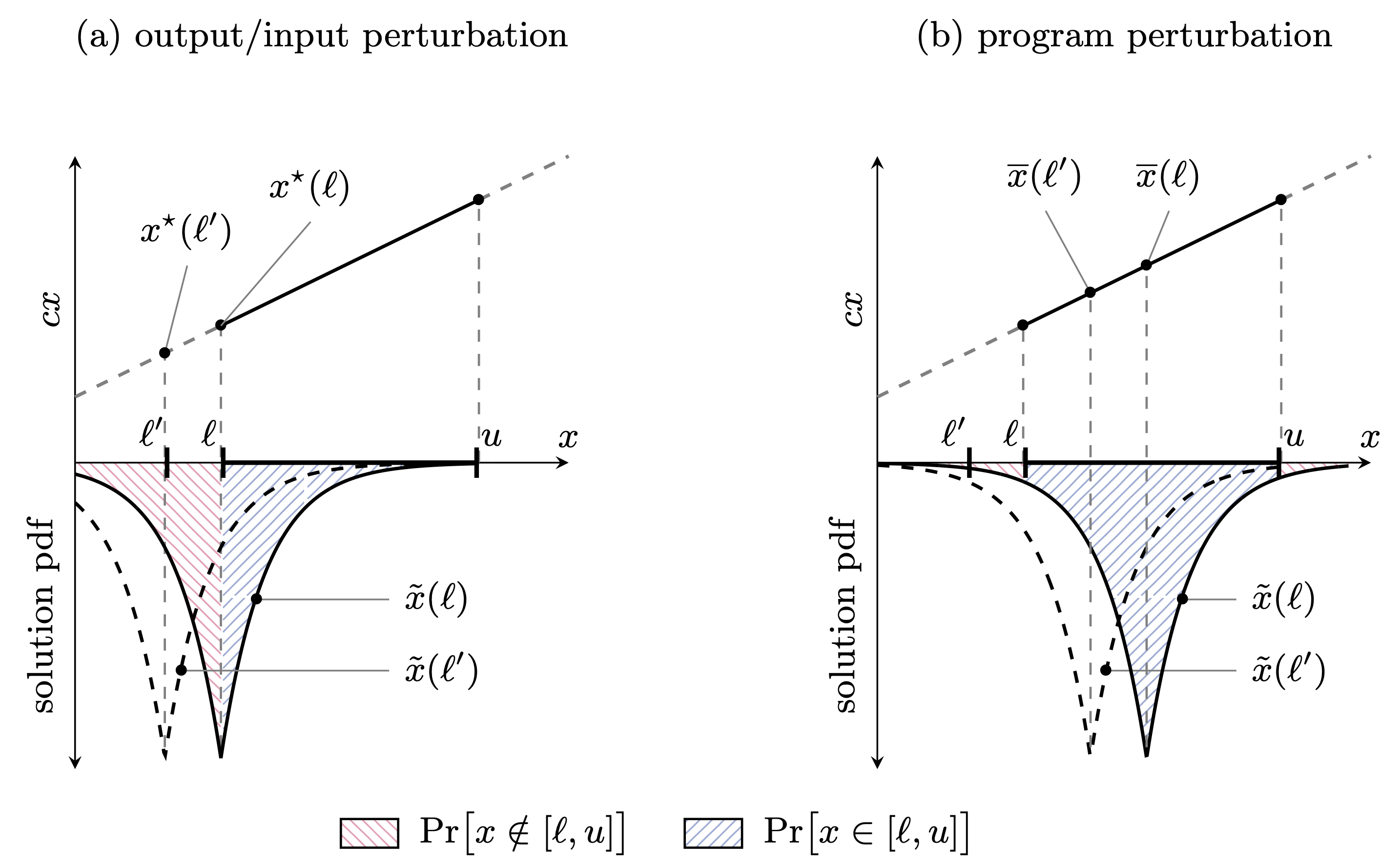}}
{
Privatization of datum $\ell$ of program \eqref{problem:LP_example} using (a) output and input perturbation strategies, and (b) program perturbation strategy. The top quadrants visualize the geometry of program \eqref{problem:LP_example}. The bottom quadrants visualize the probability density functions of perturbed optimization solutions on adjacent values $\ell$ and $\ell'$ under different perturbation strategies. \label{fig:projections}
}{}
\end{figure}

Although output and input perturbation strategies are well established in the private statistical analysis, \textcolor{maincolor}{their applications to constrained convex optimization remain limited due to feasibility requirements.} To illustrate their limits, consider a simple linear programming example. 
\begin{example}[Simple linear optimization]\normalfont\label{example:__lp__} Consider a linear optimization program
\begin{subequations}\label{problem:LP_example}
\begin{align}
    \minimize{x}\quad& cx\\
    \st\quad& \ell\leqslant x\leqslant u,
\end{align}
\end{subequations}
with one variable $x\in\mathbb{R}$ and scalar data $c,\ell,u\in\mathbb{R}_{+}$, where we assume that $\ell<u$. Let the lower limit $\ell$ be a private datum. From problem data, we see that the identity query on the optimal solution $x^{\star}=\ell$ directly exposes the private datum. The privacy goal is thus to make optimization solutions on some adjacent values $\ell$ and $\ell'$ statistically similar, such that by observing a perturbed solution, an observer would not distinguish the lower limit used in optimization. Following Proposition \ref{prop:OP_IP_dp}, this goal is achieved by perturbing either private datum $\ell$ prior to optimization or the optimal solution $x^{\star}(\ell)$ after optimization. Since $x^{\star}(\ell)=\ell$ always holds in program \eqref{problem:LP_example}, the results of the two privacy-preserving strategies are equivalent and depicted in Figure \ref{fig:projections}(a): the probability density functions (pdf's) of optimization results on adjacent values, denoted by $\tilde{x}(\ell)$ and $\tilde{x}(\ell')$, are similar up to some prescribed privacy parameters, hence achieving the privacy goal. However, the probability $\text{Pr}[x\notin[\ell,u]]$ of returning an infeasible sample from $\tilde{x}(\ell)$ is as much as 50\%. Unfortunately, the re-sampling in expectation to receive a feasible sample will degrade privacy \citep{Kifer_2011}, hence limiting the applications of output and input perturbation strategies to convex optimization.
\end{example}

Motivated by this example, our goal is to provide privacy and feasibility guarantees simultaneously by finding an alternative solution whose perturbation would be feasible. Consider, for example,  a feasible solution $\overline{x}(\ell)$ depicted in Figure \ref{fig:projections}(b). Perturbing $\overline{x}(\ell)$ with the same noise as in output and input perturbation strategies, the pdf $\tilde{x}(\ell)$ is such that the probability $\text{Pr}[x\notin[\ell,u]]$ is small. The pdf's $\tilde{x}(\ell)$ and $\tilde{x}(\ell')$ on adjacent values are also as similar as in Figure \ref{fig:projections}(a), thus achieving the privacy goal. These features, however, come at the expense of random optimality loss w.r.t. the non-perturbed, deterministic solution. To find solution $\overline{x}(\ell)$, program \eqref{problem:LP_example} must be replaced with the one, for which $\overline{x}(\ell)$ is the optimal solution given some privacy parameters and feasibility requirements. We thus refer to this strategy as program perturbation. 

We can now state the problem. We seek a random counterpart of the original, non-private convex optimization map, which satisfies the following properties:
\begin{itemize}
    \item[(a)] It produces a random optimization solution which is feasible with a high probability.
    \item[(b)] It satisfies differential privacy on adjacent optimization datasets.
    \item[(c)] Its suboptimality is controlled w.r.t. the non-private map up to prescribed parameters.
\end{itemize}

In the next section, we introduce the program perturbation strategy satisfying these properties.

\section{Differentially Private Perturbation of Convex Optimization}\label{sec:main}

\subsection{Random Optimization Map and its Chance-Constrained Optimization}\label{subsec:randoptmap}
To provide differential privacy guarantees to datasets, we introduce a zero-mean random perturbation $\bs{\zeta}\in\mathbb{R}^{k}$, and then condition the optimization solution on the realizations of $\bs{\zeta}$ by restricting the solution to the linear decision rule of the form: 
\begin{align}
    \tilde{x}(\mathcal{D})= \overline{x}(\mathcal{D}) + X(\mathcal{D})\bs{\zeta},
    \label{eq:decision_rule}
\end{align}
where $\overline{x}\in\mathbb{R}^{n}$ and $X\in\mathbb{R}^{n\times k}$ are optimization variables. Since their optimal values depend on dataset $\mathcal{D}$, sometimes we will explicitly write them as functions of datasets. Vector $\overline{x}$ is the mean (nominal) solution, which is then adjusted by the recourse $X\bs{\zeta}$ depending on realizations of $\bs{\zeta}$. Hence, by optimizing linear decision rule \eqref{eq:decision_rule}, we obtain a random optimization map. While this approach is popular in robust and stochastic programming literature to address the uncertainty of datasets \citep{ben2004adjustable,shapiro2005complexity,kuhn2011primal}, we use it to solve the reverse problem: we intentionally randomize the solution to obfuscate deterministic datasets in optimization results. 

To optimize the random optimization map $\tilde{x}(\mathcal{D})$ in order to guarantee the privacy and feasibility simultaneously, we put forth the following chance-constrained program:
\begin{subequations}\label{problem:SP}
\begin{align}
    \minimize{\overline{x},X}\quad& \mathbb{E}[c^{\top}(\overline{x}+X\bs{\zeta})]\label{problem_SP_obj}\\
    \st\quad& 
    \text{Pr}\left[b - A(\overline{x}+X\bs{\zeta}) \in\mathcal{K}\right]\geqslant1-\eta,\label{problem_SP_cone}\\
    &X\in\mathcal{X}
\end{align}
\end{subequations}
which minimizes the expected value of the cost function, subject to the joint chance constraint, which requires the satisfaction of conic constraints at least with probability $(1-\eta)$, where $\eta\in(0,1)$ is the maximum acceptable constraint violation probability. \textcolor{maincolor}{Intuitively, the feasible region for the nominal solution $\overline{x}$ is narrower in program \eqref{problem:SP} than in the deterministic problem \eqref{problem:base}, forcing $\overline{x}$ to lie within the interior of the original feasible region. This way, once the perturbation is applied to $\overline{x}^{\star}$, the resulting solution will remain within the original feasible region with probability $1-\eta$.}

The recourse variable $X$ is constrained to lie in set $\mathcal{X}$, which makes the random component of the optimization queries independent from the dataset, as required by differential privacy. The construction of this set is specific to the optimization query of interest.

\subsubsection{Optimization Queries.}

Solving optimization problems privately, the data curator answers various, application-specific queries, which are limited here to linear queries, i.e., any linear transformation of the solution vector. For example, in machine learning applications, the curator releases model parameters, such as regression weights or SVM hyperplane parameters, using an identity query. In resource allocation problems, the operator may be asked to release the total allocation of the goods or allocation costs, using a sum or weighted sum query.  To guarantee differential privacy, the decision rule in \eqref{eq:decision_rule}  must be optimized in a way which ensures that the query's random component remains independent from the dataset. In the following, we illustrate how to construct set $\mathcal{X}$ to ensure this property.

\begin{example}[Identity query]\label{example:identity_query}\normalfont The identity query returns the entire solution vector by perturbing each vector element, so the cardinality of perturbation is $|\bs{\zeta}|=|\tilde{x}|=n$. To make its random component independent from data, $X$ must be an identity matrix, i.e., $\mathcal{X}=\{X|X=I\}$, so the optimized query returns $\tilde{x}(\mathcal{D})= \overline{x}^{\star}(\mathcal{D}) + I\hat{\zeta} = \overline{x}^{\star}(\mathcal{D}) + \hat{\zeta}$, with $\hat{\zeta}$ being some realization of $\bs{\zeta}$. 
\end{example}

\begin{example}[Sum query]\label{example:sum_query}\normalfont The query returns the perturbed sum of solution vector elements, so the cardinality of perturbation is  $|\bs{\zeta}|=1$. To make its random component data independent, the only column of $X$ must sum up to 1, i.e., $\mathcal{X}=\{X|\mathbb{1}^{\top}X=1\}$, so that the optimized query returns $\mathbb{1}^{\top}\tilde{x}(\mathcal{D})= \mathbb{1}^{\top}\overline{x}^{\star}(\mathcal{D}) + \mathbb{1}^{\top}X^{\star}(\mathcal{D})\hat{\zeta}=\mathbb{1}^{\top}\overline{x}^{\star}(\mathcal{D}) + \hat{\zeta}$, with $\hat{\zeta}$ being some realization of $\bs{\zeta}$.
\end{example}

In the identity query case, set $\mathcal{X}$ makes the optimal recourse $X$ completely independent from data, i.e., $X^{\star}\neq X(\mathcal{D})$, while the sum query allows for data-dependent elements of $X$ but requires them to sum up to 1. Moreover, the cardinality $|\bs{\zeta}|$ in the identity query is larger than that of the sum query. The sum query is less restrictive, requires accommodating less noise, and is thus less expensive in terms of the expected cost. 

\subsubsection{Tractable Reformulation of Program \eqref{problem:SP}.}\label{subsec:tract_ref}

\textcolor{maincolor}{To reformulate program \eqref{problem:SP} in  a computationally tractable way, we need to address the random perturbation in the objective function and constraints. In case of linear programming, the expected cost reduces to a deterministic expression $\mathbb{E}[c^{\top}(\overline{x}+X\bs{\zeta})]\rightarrow c^{\top}\overline{x}$, thanks to the fact that $\mathbb{E}[\bs{\zeta}]=\mathbb{0}$ by design. The expected cost of non-linear objective functions, such as the case of QP and SDP, is explained later in Section \ref{sec:application}.}


\textcolor{maincolor}{Addressing the joint chance constraint \eqref{problem_SP_cone}, we can adopt either a sample-based reformulation or, in some special cases, an analytic reformulation. In the baseline sample-based approach, first developed by \citet{calafiore2005uncertain}, the entries of \eqref{problem_SP_cone} are enforced on a finite number of samples $\hat{\zeta}_{1},\dots,\hat{\zeta}_{S}$ extracted from the distribution of $\boldsymbol{\zeta}$, resulting in the sampled convex program (SCP) whose solution (if exists) is feasible for the original chance-constrained program \eqref{problem:SP}, provided that a sufficient number $S$ of samples is drawn. Unfortunately, the sample size requirement $S$ significantly increases in the number of decision variables; for example, \cite[Remark 6]{calafiore2006scenario} states $S\gg10,000$ for a problem with 10 decision variables. As a result, a straightforward constraint sampling may produce a program approximation which can be challenging to solve by standard interior-point solvers, especially when the set $\mathcal{K}$ includes positive semidefinite cones.} 

\textcolor{maincolor}{To reduce the size of the SCP, \citet{margellos2014road} proposed enforcing the entries of the chance constraint \eqref{problem_SP_cone} at the {\it vertices} of the hyperrectangular uncertainty set built upon extracted samples $\hat{\zeta}_{1},\dots,\hat{\zeta}_{S}$. The virtue of Margellos' approach is that the required sample size $S$ does not depend on the number $n$ of decision variables but rather on the dimension $k$ of the random perturbation $\boldsymbol{\zeta}$. Thus, when working with aggregated (e.g., sum, average) optimization queries for which $k\ll n$, the resulting SCP could feature a smaller size than that obtained under \citet{calafiore2005uncertain}'s approach. Moreover, because the entries of $\boldsymbol{\zeta}$ are uncorrelated, the hyperrectangular uncertainty set is not as conservative compared to the cases with correlated random variables. The downside of Margellos' approach, however, is that extracting $2^{k}$ vertices from the hyperrectangular set could itself be a challenging task for large $k$ \citep{khachiyan2009generating}. The reformulated chance-constrained program then takes the form:}
\begin{subequations}\label{problem:SP_tractable}
\begin{align}
    \minimize{\overline{x},X\in\mathcal{X}}\quad& c^{\top}\overline{x}\label{problem_SP_tr_obj}\\
    \st\quad& 
    b - A(\overline{x}+X\zeta_{i}^{v}) \in\mathcal{K},\quad\forall i=1,\dots,2^{k},\label{problem_SP_tr_cone}
\end{align} 
\end{subequations}
where $\zeta_{i}^{v}\in\mathbb{R}^{k}$ denotes a particular vertex $i$ \textcolor{maincolor}{of the hyperrectangular set. With a carefully chosen sample size $S$, constraint \eqref{problem_SP_tr_cone} approximates \eqref{problem_SP_cone} with some confidence level} as per the following result from \citep{margellos2014road}.
\begin{proposition}[Sample size requirement 1]\label{prop:sample_size_1}\normalfont 
 \textcolor{maincolor}{Let $\beta\in(0,1)$ be some confidence parameter, and let $\zeta_{1}^{v},\dots,\zeta_{2^{k}}^{v}$ be the vertices of the hyperrectangular uncertainty set constructed upon $S$ extracted samples. If the sample size is chosen according to  
$$
    S \geqslant \left\lceil
    \frac{1}{\eta}\frac{e}{e-1}\left(2k-1+\text{ln}\frac{1}{\beta}\right)
    \right\rceil,
$$
where $e$ is the Euler’s number and operator $\lceil x\rceil$ returns the least integer greater than $x$, then the optimal solution to problem \eqref{problem:SP_tractable}, if exists, is a feasible solution to the chance-constrained program \eqref{problem:SP} with probability at least $(1-\beta)$.}
\end{proposition}

This reformulation is universal in accommodating any constraint type in the joint chance constraint, including linear, quadratic, semidefinite, exponential and other constraints. It also preserves the problem type, e.g., the linear program remains linear when the perturbation is applied. In some special cases below, however, it makes sense to use an alternative, analytic reformulation. 
\begin{remark}[Linear constraints]\label{remark:lin_con_cc}\normalfont For large $k$, as in the identity query case, vertex enumeration becomes computationally challenging \citep{khachiyan2009generating}. However, in the case of linear constraints in \eqref{problem_SP_cone}, i.e., $\text{Pr}[b - A(\overline{x}+X\bs{\zeta}) \in \mathbb{R}_{\textcolor{maincolor}{+}}^{m}]\geqslant1-\eta$, we can reformulate the joint chance constraint analytically into a set of efficient second-order cone constraints. \textcolor{maincolor}{Specifically, each entry of the joint chance constraint is treated as an individual chance constraint, i.e., 
\begin{align}\label{eq:ind_cc}
\text{Pr}[b_{i} - A_{i}^{\top}(\overline{x}+X\bs{\zeta}) \in \mathbb{R}_{+}]\geqslant1-\overline{\eta}_{i},\quad\forall i=1,\dots,m,
\end{align}
where $\overline{\eta}_{1},\dots,\overline{\eta}_{m}\in(0.5,1)$ are individual constraint violation probabilities. When they satisfy $\sum_{i=1}^{m}\overline{\eta}_{i}=\eta$, the set of individual constraints in \eqref{eq:soc_ref} guarantees the satisfaction of the original joint chance constraint, as per Bonferroni approximation \citep{xie2019optimized}. Then, each individual chance-constraint can be reformulated as a second-order cone constraint of the form 
\begin{align}\label{eq:soc_ref}
z(\overline{\eta}_{i})\norm{F[AX]_{i}^{\top}}\leqslant b_{i} - A_{i}^{\top}x
\quad\forall i=1,\dots,m,
\end{align} where $F$ is the factorization of the \textcolor{maincolor}{covariance matrix $\Sigma$ of $\boldsymbol{\zeta}$}, i.e., $\Sigma=FF^{\top}$, and $z(\cdot)$ is the ``safety'' parameter in terminology of \citet{nemirovski2007convex}. For any distribution, e.g., Laplace distribution, setting $z(\overline{\eta}_{i})=\sqrt{(1-\overline{\eta}_{i})/\overline{\eta}_{i}}$ in \eqref{eq:soc_ref} conservatively approximates \eqref{eq:ind_cc}, as per Chebyshev inequality \citep{xie2017distributionally}. For the special case of Gaussian perturbation, reformulation \eqref{eq:soc_ref} is exact when $z(\overline{\eta}_{i})$ amounts to the inverse cumulative distribution function of the standard Gaussian distribution at $(1-\overline{\eta}_{i})^{\text{th}}$ quantile.}
\end{remark}
\begin{remark}[Equality constraint]\label{remark:eq_con}\normalfont The perturbation of equality constraints in the linear decision rule takes the form  $b - A(\overline{x}+X\bs{\zeta}) = \mathbb{0}$. The perturbed equality is satisfied with probability 1 when replaced with two deterministic equations $b-A\overline{x}= \mathbb{0}$ and $AX=\mathbb{0}$.
\end{remark}

\textcolor{maincolor}{For a comprehensive reviewer of alternative techniques to reformulate program \eqref{problem:SP} in a computationally tractable way, the interested reader is kindly referred to \cite{geng2019data} and \cite{kuccukyavuz2021chance}.}

\subsection{Differential Privacy Guarantees}\label{subsec:privacy_theorems}

\textcolor{maincolor}{Here, we establish formal privacy guarantees. For the cases where the sensitivity of the solution map is known or can be upper-bounded by domain knowledge, we provide the pure $\varepsilon-$DP guarantee. For the cases when the exact sensitivity is unknown or hard to upper bound, we provide the sample-based estimation method and the corresponding probabilistic DP guarantee. 
}

\subsubsection{Pure Differential Privacy Guarantee:}

\textcolor{maincolor}{By calibrating the Laplace perturbation to the worst-case sensitivity of optimization queries, we can establish the following differential privacy guarantees. The former is provided for the identity query which privately releases the entire solution vector $x$, and the latter is for the sum queries typically used for aggregated statistics.}

\begin{theorem}[$\boldsymbol{\varepsilon-}$DP identity optimization query]\label{th:eps_dp_id}\normalfont Let $\bs{\zeta}\sim\text{Lap}(\Delta_{1}/\varepsilon)^{n}$, where $\Delta_{1}$ is the worst-case $\ell_{1}-$sensitivity of the identity query on $\alpha-$adjacent optimization datasets $\mathcal{D},\mathcal{D}'\in\mathbb{D}$. Then, if program \eqref{problem:SP} returns the optimal solution \textcolor{maincolor}{($\bar{x}^\star, X^\star$)}, identity query $\overline{x}^{\star}(\mathcal{D}) + X^{\star}(\mathcal{D})\hat{\zeta}$ with $\hat{\zeta}$ being some realization of $\bs{\zeta}$, is $\varepsilon-$differentially private. That is, for all $\hat{x}$ in the query's range:
\begin{align*}
    \text{Pr}[\overline{x}^{\star}(\mathcal{D}) + X(\mathcal{D})^{\star}\bs{\zeta} \in \hat{x}]\leqslant\text{Pr}[\overline{x}^{\star}(\mathcal{D}') + X^{\star}(\mathcal{D}')\bs{\zeta} \in \hat{x}]\text{exp}(\varepsilon). 
\end{align*}
\end{theorem}

\begin{theorem}[$\boldsymbol{\varepsilon-}$DP sum optimization query]\label{th:eps_dp_sum}\normalfont Let $\bs{\zeta}\sim\text{Lap}(\Delta_{1}/\varepsilon)^{1}$, where $\Delta_{1}$ is the worst-case $\ell_{1}-$sensitivity of the sum query on $\alpha-$adjacent optimization datasets $\mathcal{D},\mathcal{D}'\in\mathbb{D}$. Then, if program \eqref{problem:SP} returns the optimal solution \textcolor{maincolor}{($\bar{x}^\star, X^\star$)}, sum query $\mathbb{1}^{\top}(\overline{x}^{\star}(\mathcal{D}) + X^{\star}(\mathcal{D})\hat{\zeta})$ with $\hat{\zeta}$ being some realization of $\bs{\zeta}$, is $\varepsilon-$differentially private. That is, for all $\hat{x}$ in the query's range:
\begin{align*}
    \text{Pr}[\mathbb{1}^{\top}(\overline{x}^{\star}(\mathcal{D}) + X(\mathcal{D})^{\star}\bs{\zeta}) \in \hat{x}]\leqslant\text{Pr}[\mathbb{1}^{\top}(\overline{x}^{\star}(\mathcal{D}') + X^{\star}(\mathcal{D}')\bs{\zeta}) \in \hat{x}]\text{exp}(\varepsilon). 
\end{align*}
\end{theorem}

\textcolor{maincolor}{Interestingly, the private sum query can also be answered by summing the outcome of the private identity query, as per Theorem \ref{th:post_processing} of post-processing immunity. However, this would be sub-optimal as the identity query requires substantially more noise then the sum query. Indeed, the identity query requires $n-$dimensional perturbation as opposed to $1-$dimensional perturbation of the sum query; moreover, the sensitivity of identity queries are larger than those of aggregated queries.}

\subsubsection{Probabilistic Differential Privacy Guarantee:}

\textcolor{maincolor}{Often, the worst-case sensitivity $\Delta_{p}$ of the solution map $x(\mathcal{D})$ to adjacent datasets is unknown and requires estimation, which} is itself a non-convex, hierarchical optimization problem, e.g., for the identity query it takes the form:
\begin{subequations}\label{prog:sens_non_convex}
\begin{align}
    \Delta_{p}=\maximize{\mathcal{D},\mathcal{D}'\in\mathbb{D}}\quad&\norm{x(\mathcal{D}) - x(\mathcal{D}')}_{p}\\
    \st\quad&\norm{\mathcal{D}-\mathcal{D}'}\leqslant\alpha,
\end{align}
\end{subequations}
where $x(\mathcal{D})$ and $x(\mathcal{D}')$ are two embedded optimization programs. For this problem, there is no efficient algorithm to find the global optimal solution. However, following the approach of \cite{dobbe2018customized}, we can frame this problem as the following chance-constrained stochastic program:
\begin{subequations}\label{prog:sens_cc_program}
\begin{align}
    \minimize{t}\quad&t\\
    \st\quad&\text{Pr}_{\textcolor{maincolor}{\mathbb{P}_{\mathcal{D}\sim_{\alpha}\mathcal{D}'}}}\left[\norm{x(\mathcal{D}) - x(\mathcal{D}')}_{p} \leqslant t\right]\geqslant 1-\gamma, 
\end{align}
\end{subequations}
in one variable $t$, where the probability in the chance constraint is taken w.r.t. the joint distribution \textcolor{maincolor}{$\mathbb{P}_{\mathcal{D}\sim_{\alpha}\mathcal{D}'}$ of $\alpha-$adjacent datasets $\mathcal{D}$ and $\mathcal{D}'$}. Solving program \eqref{prog:sens_cc_program} results in solution sensitivity which is valid for $1-\gamma$ instances of $\alpha-$adjacent dataset pairs $(\mathcal{D},\mathcal{D}')$, leaving out a small percentage $\gamma$ of dataset pairs on which the sensitivity could exceed solution $t^{\star}$; \textcolor{maincolor}{hence, $t^{\star}$ is the lower bound on $\Delta_{p}$}. We can solve this problem efficiently by optimizing variable $t$ on a finite number $S$ of adjacent datasets \textcolor{maincolor}{from distribution $\mathbb{P}_{\mathcal{D}\sim_{\alpha}\mathcal{D}'}$}, i.e., using the following program:
\begin{subequations}\label{prog:sens_sample}
\begin{align}
    \minimize{t}\quad&t\\
    \st\quad&\norm{x(\mathcal{D}_{s}) - x(\mathcal{D}_{s}')}_{p} \leqslant t, \quad \forall s=1,\dots,S,
\end{align}
\end{subequations}
where $(\mathcal{D}_{s},\mathcal{D}_{s}')$ is sampled pair of adjacent datasets and $x(\cdot)$ is the solution to program \eqref{problem:base} evaluated offline, and $S$ is the sample size. The optimal solution $t^{\star}$ is the lower bound on the worst-case sensitivity, which can be tighten by increasing the sample size. \textcolor{maincolor}{Formulation \eqref{prog:sens_sample} is the sampled convex program in the sense of \cite[Definition 1]{calafiore2005uncertain} and amendable to Corollary 1 from \citet{calafiore2005uncertain} on the minimum sample size requirement that helps us to statistically lower bound the worst-case solution sensitivity as follows.}

\begin{proposition}[Sample size requirement 2]\label{prop:sample_size}\normalfont \textcolor{maincolor}{Given probability parameter $\gamma\in(0,1)$ and confidence parameter $\beta\in(0,1)$, if the sample size $S\geqslant 1/(\gamma\beta)-1$, the optimal solution $t^{\star}$ to problem  \eqref{prog:sens_sample} is a feasible solution to the chance-constrained program \eqref{prog:sens_cc_program} with probability at least $(1-\beta)$.}
\end{proposition}

This requirement is independent from the joint distribution of adjacent datasets, i.e., particularly from the dimension of the joint distribution, which makes the sensitivity computation tractable. \textcolor{maincolor}{
Since $t^{\star}$ is the lower bound on the true sensitivity, the perturbation calibrated to $t^{\star}$ will not satisfy the pure $\varepsilon-$DP, but its relaxed, probabilistic version, termed $(\varepsilon,\gamma)-$probabilistic differential privacy ($\varepsilon,\gamma-$PDP) in the sense that the $\varepsilon-$DP is provided for $(1-\gamma)\%$ of $\alpha-$adjacent datasets as per the following result.}

\textcolor{maincolor}{
\begin{theorem}[$\boldsymbol{(\varepsilon,\gamma)-}$PDP identity optimization query]\label{th:eps_pdp_id}\normalfont Let $\bs{\zeta}\sim\text{Lap}(t^{\star}/\varepsilon)^{n}$, where $t^{\star}$ is the solution to the chance-constrained program \eqref{prog:sens_cc_program} for $p=1$. Then, if program \eqref{problem:SP} returns the optimal solution \textcolor{maincolor}{($\bar{x}^\star, X^\star$)}, the identity query $\overline{x}^{\star}(\mathcal{D}) + X^{\star}(\mathcal{D})\hat{\zeta}$ with $\hat{\zeta}$ being some realization of $\bs{\zeta}$, satisfies probabilistic $\varepsilon-$differential privacy in the sense that the following inequality holds 
\begin{align*}
    \text{Pr}[\overline{x}^{\star}(\mathcal{D}) + X(\mathcal{D})^{\star}\bs{\zeta} \in \hat{x}]\leqslant\text{Pr}[\overline{x}^{\star}(\mathcal{D}') + X^{\star}(\mathcal{D}')\bs{\zeta} \in \hat{x}]\text{exp}(\varepsilon). 
\end{align*}
\end{theorem}
for $(1-\gamma)\%$ of $\alpha-$adjacent dataset pairs $(\mathcal{D},\mathcal{D}')\in\mathbb{D}.$}

Similarly, we can provide the $(\varepsilon,\gamma)-$PDP guarantee for the sum query. 
\textcolor{maincolor}{
\begin{remark}
    Theorem \ref{th:eps_pdp_id} builds on the optimal solution of the stochastic program \eqref{prog:sens_cc_program}. In practice, we would solve a sample approximation problem \eqref{prog:sens_sample} instead. As per Proposition \ref{prop:sample_size}, the sample size $S$ depends on probability parameter $\gamma$ and confidence parameter $\beta$, both affecting $S$, and thus the sensitivity estimation. Therefore, due to sampling nature of the sensitivity estimation, we also need to specify the level of confidence when establishing the PDP guarantee. 
\end{remark}}

\subsubsection{Approximate Differential Privacy Guarantee}




The pure DP guarantee can be relaxed to $(\varepsilon,\delta)$-differential privacy by calibrating the perturbation to the $\ell_{2}-$sensitivities and choosing a Gaussian distribution instead of Laplace, similar to the Gaussian mechanism in \citep[Theorem A.1]{dwork2014algorithmic}. \textcolor{maincolor}{Although the variance of the Gaussian mechanism is larger than that of the Laplace mechanism, the analytic reformulation of individual chance constraints (see Remark \ref{remark:lin_con_cc}) is exact for the Gaussian perturbation, potentially resulting in a less conservative approximation of the chance-constrained program \eqref{problem:SP}.}
\begin{theorem}[$\boldsymbol{(\varepsilon,\delta)-}$differentially private identity query]\label{th:eps_delta_dp_id}\normalfont 
Let $\bs{\zeta}\sim N(\mathbb{0},\text{diag}[\mathbb{1}\cdot\sigma^2])$, where $\sigma=\sqrt{2\text{ln}(1.25/\delta)}\Delta_{2}/\varepsilon$, where $\Delta_{2}$ is the  $\ell_{2}-$sensitivity of identity query to adjacent datasets $\mathcal{D},\mathcal{D}'\in\mathbb{D}$. Then, if program \eqref{problem:SP} returns the optimal solution $\bar{x}^*, X^*$, identity query $\overline{x}^{\star}(\mathcal{D}) + X^{\star}(\mathcal{D})\hat{\zeta}$ with $\hat{\zeta}$ being some realization of $\bs{\zeta}$, is $(\varepsilon,\delta)-$differentially private. That is, for all $\hat{x}$ in the query's range:
\begin{align*}
    \text{Pr}[\overline{x}^{\star}(\mathcal{D}) + X(\mathcal{D})^{\star}\bs{\zeta} \in \hat{x}]\leqslant\text{Pr}[\overline{x}^{\star}(\mathcal{D}') + X^{\star}(\mathcal{D}')\bs{\zeta} \in \hat{x}]\text{exp}(\varepsilon) + \delta.
\end{align*}
\end{theorem}
Similarly, the Gaussian mechanism is built for the privacy-preserving sum query.

\subsection{Controlling Suboptimality of the Random Optimization Map}\label{subsec:CVaR}

Applying privacy-preserving perturbations inevitably results in the optimality loss w.r.t. the non-private solution, which itself is a random variable. For some generic optimization objective function $\ell:\mathbb{R}^{n}\mapsto\mathbb{R}$, the expected optimality loss amounts to
\begin{align}
    \mathbb{E}\big[\norm{\ell(\overline{x}^{\star}+X^{\star}\bs{\zeta}) - \ell(x^{\star})}\big] = \mathbb{E}\big[\ell(\overline{x}^{\star}+X^{\star}\bs{\zeta}) - \ell(x^{\star}) \big]=\mathbb{E}\big[\ell(\overline{x}^{\star}+X^{\star}\bs{\zeta}) \big] - \ell(x^{\star}),
\end{align} 
i.e., to the distance between optimized deterministic \eqref{problem:base} and random  \eqref{eq:decision_rule} solution maps. Therefore, optimizing the solution map using program \eqref{problem:SP} guarantees the minimal optimality loss in expectation. However, in many applications of interest, any optimality loss directly translates to financial losses, which decision-makers would prefer to hedge against using alternative measures of risks, such as value-at-risk (VaR) or conditional VaR (CVaR). In our notation, VaR is defined as
\begin{subequations}
\begin{align}
    \text{VaR}_{1-q}\left[\ell(\overline{x}^{\star}+X^{\star}\bs{\zeta})\right] = \text{inf}\{\gamma\in\mathbb{R}\;|\;\text{Pr}[\ell(\overline{x}^{\star}+X^{\star}\bs{\zeta})\leqslant\gamma]\geqslant1-q\},
\end{align}
i.e., as the minimal value $\gamma^{\star}$ that objective function $\ell$ will not exceed with probability $1-q$. While VaR's interpretation is appealing, it is not a coherent risk measure and thus hard to optimize \citep{rockafellar2000optimization}. We hence optimize the CVaR, which in our notation amounts to 
\begin{align}
    \text{CVaR}_{1-q}\left[\ell(\overline{x}^{\star}+X^{\star}\bs{\zeta})\right] = \frac{1}{1-q}\int_{1-q}^{1}\text{VaR}_{\beta}\left[\ell(\overline{x}^{\star}+X^{\star}\bs{\zeta})\right]d\beta,\label{eq:CVaR}
\end{align}
i.e., the average $(1-q)-$quantile of the random objective value. It is shown by \citet{rockafellar2000optimization} that \textcolor{maincolor}{CVaR \eqref{eq:CVaR} can approximately be optimized by minimizing its sample approximation, namely}
\begin{align}\label{cvar:sample_approx_}
    \minimize{\gamma}\quad&
    \gamma + \frac{1}{(1-q)S}\sum_{s=1}^{S}\big[\ell(\overline{x}^{\star}+X^{\star}\hat{\zeta}_{s}) - \gamma\big]^{+},
\end{align}
\end{subequations}
where $\hat{\zeta}_{1},\dots,\hat{\zeta}_{S}$ being a set of samples extracted from the distribution of $\bs{\zeta}$, and $[x]^{+}$ is the projection of $x$ on the non-negative space. This gives rise to the following convex CVaR optimization using a set of auxiliary variables $z_{1},\dots,z_{S}$:
\begin{subequations}\label{prog:cvar}
\begin{align}
    \minimize{\gamma,z_{1},\dots,z_{S}}\quad&
    \gamma + \frac{1}{(1-q)S}\sum_{s=1}^{S}z_{s}\\
    \st\quad&z_{s}\geqslant0,\quad z_{s}\geqslant\ell(\overline{x}^{\star}+X^{\star}\hat{\zeta}_{s}) - \gamma,\quad\forall s=1,\dots,S.
\end{align}
\end{subequations}
Therefore, co-optimizing the linear decision rule and CVaR variables yields the optimized solution map which is optimal w.r.t. the average over $(1-q)\%$ worst-case scenarios of optimality loss. This distribution-agnostic CVaR optimization is suitable to accommodate both Laplacian and Gaussian perturbations of differential privacy. Moreover, it does not violate the privacy guarantees of Theorems \ref{th:eps_dp_id}--\ref{th:eps_delta_dp_id} as long as the query-specific constraint $\mathcal{X}$ is enforced on recourse variable $X$.

\textcolor{maincolor}{\subsection{On the Feasibility of the Program Perturbation Strategy}}

\textcolor{maincolor}{In the proposed privacy-preserving framework, the curator specifies the optimization query, calibrates the noise w.r.t. the desired privacy requirements, and then solves an approximation (e.g., sampled or analytic) of the chance-constrained program \eqref{problem:SP}. The privacy guarantees established in Theorems \ref{th:eps_dp_id} through \ref{th:eps_delta_dp_id} depend on whether the approximation of the chance-constrained program \eqref{problem:SP} returns a feasible solution. If the approximation is feasible, the curator samples and provides a private response to the query. However, if the approximation is infeasible, then either (a) the chosen approximation of the chance constraint is overly conservative, e.g., the sampled program \eqref{problem:SP_tractable} is infeasible, but the original chance-constrained program  \eqref{problem:SP} has a solution; or (b) the original chance-constrained program itself does not have a solution. In case (a), the curator fails to specify an approximation that allows for the retrieval of the privacy-preserving solution. In case (b), the required level of privacy preservation is not achievable, i.e., there is excessive noise for the limited feasible region. In either cases, the infeasibility status confirms the impossibility of privacy preservation with the specified parameters.}

\textcolor{maincolor}{However, in some special cases, the possibility of accommodating a privacy-preserving perturbation within the feasible region can be assessed {\it a priori}. For instance, consider the linear program in Example \ref{example:__lp__} with a double-sided constraint $\ell \leqslant x \leqslant u$ on a scalar variable $x$. This can be done by relating the diameter $d_{\text{LP}}=u-\ell$ of the feasible region with the length of the interval containing $(1-\eta)\%$ of the probability mass of the Laplace random variable. For example, for $\eta=5\%$, the interval containing 95\% of the probability mass has a length of $d_{\text{lap}} = -2b\ln(0.05)$ -- the expression derived from solving the CDF of the Laplace distribution at point $x$, with $b$ being the scale of the Laplace distribution. Then, after setting $b=\frac{\Delta_{1}}{\alpha}$, if $d_{\text{lap}} \leqslant d_{\text{LP}}$, the chance-constrained problem has a feasible solution. The challenge then becomes finding a good approximation of the chance-constrained program, e.g., a sampled convex program, that can retrieve it. For polyhedral constraint set, the interval containing $(1-\eta)\%$ of the probability mass can be related to the smallest radius of the maximum-volume ellipsoid inscribed in the polyhedral set \citep{ben1999robust}. However, such {\it a priori} certification remains very challenging for arbitrary convex programs.}

\vspace{0.5cm}
\section{Applications}\label{sec:application}

We demonstrate a privacy-preserving perturbation of linear optimization programs using the optimal power flow problem in Section \ref{subsec:OPF}. We also apply the framework to classification and regression problems in energy systems to illustrate differentially private quadratic programming in Sections \ref{subsec:SVM} and \ref{subsec:Regression}. Last, we present a private semidefinite programming using an ellipsoid fitting problem in Section \ref{subsec:elliposid}. We present programs in their explicit form, but all of them are representable through the reference program \eqref{problem:base} as shown in e-companion \ref{app:conic_representations}. For datasets, codes and implementation details we also refer to the online repository associated with this paper \citep{dvorkin2022github}. 

\subsection{Private Linear Programming: DC Optimal Power Flow (OPF) Problem}\label{subsec:OPF}

Operation of an electric power network with $N$ nodes and $E$ transmission lines is governed by the solution of the optimal power flow problem, which identifies the least-cost generation allocation $x\in\mathbb{R}^{N}$ that satisfies electricity demand $d\in\mathbb{R}^{N}$ and technical limits on generation $x^{\text{min}},x^{\text{max}}\in\mathbb{R}^{N}$ and power flows $f^{\text{max}}\in\mathbb{R}^{E}$. Using the DC power flow representation, the nodal power mismatch $x-d$ can be translated to the power flows using a matrix $F\in\mathbb{R}^{E\times N}$ of power transfer distribution factors \citep{chatzivasileiadis2018lecture}. Then, the problem takes the following linear programming form:
\begin{subequations}\label{OPF:det}
\begin{align}
    \minimize{x}\quad&c^{\top}x\label{OPF:obj}\\
    \st\quad&\mathbb{1}^{\top}(x-d)=0,\label{OPF:det_bal}\\
    &|F(x-d)|\leqslant f^{\text{max}},\label{OPF:det_flow}\\
    &x^{\text{min}}\leqslant x\leqslant x^{\text{max}},\label{OPF:det_gen}
\end{align}
\end{subequations}
where the generation cost function is minimized subject to power balance constraint \eqref{OPF:det_bal} and limits on power flows and generation in \eqref{OPF:det_flow} and \eqref{OPF:det_gen}, respectively. 

The privacy goal here is to protect any individual element of demand vector $d$ up to some adjacency value $\alpha$ when answering queries on the optimal objective function value. \textcolor{maincolor}{For example, here $\alpha=1$ MWh means that any two vectors $d$ and $d'$ are different in one individual electricity demand by at most $1$ MWh.} 
In this case, \textcolor{maincolor}{the sensitivity of the optimal objective function query is bounded by} $\textcolor{maincolor}{c}_{\textcolor{maincolor}{\text{max}}}\alpha$, i.e., by the change in generation cost of the most expensive generator with cost $c_{\text{max}}=\text{max}\{c_{1},\dots,c_{N}\}$. \textcolor{maincolor}{From Proposition \ref{prop:OP_IP_dp}, the output perturbation strategy then takes the form $c^{\top}x^{\star} + \text{Lap}(c_{\text{max}}/\varepsilon)$}. \textcolor{maincolor}{The input perturbation strategy, instead, would first obfuscate the load vector before optimization. The sensitivity of the private identity query when adjacent vectors are different by at most one element is bounded by $\alpha$, and the input perturbation takes the form: load obfuscation $\tilde{d} = d + \text{Lap}(\alpha/\varepsilon)^{N}$, followed by optimization of \eqref{OPF:det} on $\tilde{d}$.} For the output perturbation, however, there is no formal guarantee that the perturbed objective function value is feasible, i.e., that there exists a feasible allocation $x$ that provides the perturbed objective function value. For the input perturbation, there is also no guarantee that the perturb load vector admits a feasible solution. To provide feasibility and privacy guarantees simultaneously, our proposed program perturbation strategy requires solving the following chance-constrained counterpart of program \eqref{OPF:det}:
\begin{subequations}\label{OPF:sto}
\begin{align}
    \minimize{\overline{x},X}\quad&\mathbb{E}[c^{\top}(\overline{x}+X\bs{\zeta})]\\
    \st\quad&\mathbb{1}^{\top}(\overline{x}+X\bs{\zeta}-d)=0,\label{sto_opf:eq_con_infty}\\
    &
    \text{Pr}\left[
    \begin{aligned}
    &|F(\overline{x}+X\bs{\zeta}-d)|\leqslant f^{\text{max}}\\
    &x^{\text{min}}\leqslant \overline{x}+X\bs{\zeta}\leqslant x^{\text{max}}
    \end{aligned}
    \right]
    \geqslant1-\eta,\label{OPF_sto_cc}\\
    &c^{\top}X=1,\label{OPF_sto_recourse}
\end{align}
\end{subequations}
in variables $\overline{x}$ and $X$, where $\bs{\zeta}\sim\text{Lap}(\textcolor{maincolor}{c_{\text{max}}}/\varepsilon)$. \textcolor{maincolor}{Towards its tractable reformulation, recall that $\mathbb{E}[\boldsymbol{\zeta}]=0$, so that the objective function reduces to $c^{\top}\overline{x}$. Following Remark \ref{remark:eq_con}, the stochastic equality constraint \eqref{sto_opf:eq_con_infty} is substituted with an equivalent set of two deterministic constraints: $\mathbb{1}^{\top}(\overline{x}-d)=0$ and $\mathbb{1}^{\top}X=\mathbb{0}$. The joint chance constraint \eqref{OPF_sto_cc} is addressed by the sample approximation approach outlined in Section \ref{subsec:tract_ref}. To provide a feasibility guarantee, the maximum theoretical constraint violation probability, $\eta$, is set to $1\%$, and we select the confidence parameter $\beta$ at $10\%$. Following Proposition \ref{prop:sample_size_1}, the vertices of the $k-$dimensional uncertainty set, with $k=1$, were computed using the sample size requirement $S=523$.} The last constraint \eqref{OPF_sto_recourse} is enforced to require data-independence of the objective function's recourse to enable differential privacy guarantees, similar to Example \ref{example:sum_query}. Once program \eqref{OPF:sto} is solved to optimality, the differentially private objective function release takes the form $c^{\top}\overline{x}^{\star} \textcolor{maincolor}{+\text{Lap}(c_{\text{max}}/\varepsilon)}$.

\begin{table}
\TABLE
{Summary of the three differential privacy strategies for private OPF cost queries. Optimality loss and query infeasibility (inf) are given in percentage and evaluated on 1000 realizations of the random perturbation.\label{tab:OPF_summary}}
{
\footnotesize
\begin{tabular}{llcccccc}
\toprule
\multirow{3}{*}{dataset} & \multicolumn{1}{c}{\multirow{3}{*}{\begin{tabular}[l]{@{}l@{}}perturbation \\ strategy\end{tabular}}} & \multicolumn{6}{c}{dataset adjacency (MWh)} \\
\cmidrule(lr){3-8}
 & \multicolumn{1}{c}{} & \multicolumn{2}{c}{
 \begin{tabular}[c]{@{}c@{}}$\alpha=1$ \\ {\footnotesize (low-privacy)} \end{tabular}
 } & \multicolumn{2}{c}{
\begin{tabular}[c]{@{}c@{}}$\alpha=3$ \\ {\footnotesize (mid-privacy)} \end{tabular}
 } & \multicolumn{2}{c}{
\begin{tabular}[c]{@{}c@{}}$\alpha=10$ \\ {\footnotesize (high-privacy)} \end{tabular}
 } \\
\cmidrule(lr){3-4}\cmidrule(lr){5-6}\cmidrule(lr){7-8}
 & \multicolumn{1}{c}{} & loss & inf & loss & inf & loss & inf \\
\midrule
\multirow{3}{*}{5\_pjm} & input & 0.00 & 48.90 & 0.10 & 51.20 & 0.10 & 50.20 \\
 & output & 0.00 & 51.20 & 0.00 & 49.60 & 0.10 & 51.60 \\
 & program & 1.07 & 0.50 & 7.00 & 0.00 & 12.10 & 0.20 \\
\midrule
\multirow{3}{*}{14\_ieee} & input & 0.00 & 49.00 & 0.30 & 49.00 & 0.40 & 52.90 \\
 & output & 0.10 & 48.50 & 0.20 & 48.10 & 0.60 & 52.00 \\
 & program & 7.10 & 0.40 & 25.20 & 0.30 & -- & -- \\
\midrule
\multirow{3}{*}{24\_ieee} & input & 0.00 & 51.50 & 0.00 & 49.90 & 0.20 & 50.30 \\
 & output & 0.10 & 52.70 & 0.00 & 51.50 & 0.10 & 48.80 \\
 & program & 1.70 & 0.10 & 5.10 & 0.00 & 17.10 & 0.00 \\
\midrule
\multirow{3}{*}{57\_ieee} & input & 0.00 & 50.20 & 0.00 & 51.10 & 0.20 & 55.90 \\
 & output & 0.00 & 49.40 & 0.00 & 49.10 & 0.05 & 49.80 \\
 & program & 0.70 & 0.10 & 2.20 & 0.10 & 6.70 & 0.20 \\
\midrule
\multirow{3}{*}{89\_pegase} & input & 0.00 & 50.60 & 0.00 & 55.00 & 0.30 & 67.50 \\
 & output & 0.00 & 48.50 & 0.00 & 51.10 & 0.00 & 48.30 \\
 & program & 0.30 & 0.00 & 0.80 & 0.10 & 2.50 & 0.10 \\
\bottomrule
\end{tabular}
}
{}
\end{table}

For numerical experiments, we use a series of benchmark power networks from \citep{coffrin2018powermodels} and choose to provide $1-$differential privacy ($\varepsilon=1$) for varying adjacency parameter $\alpha$. The input, output and program perturbation strategies are compared in terms of the expected suboptimality gap (w.r.t. deterministic, non-private solution) and probability of infeasible query in Table \ref{tab:OPF_summary}. The results show a high infeasibility of input and output perturbation strategies: the former does not necessarily produce a feasible optimization dataset, and the latter will always fail to produce a feasible query with the value below the optimal one, i.e., for negative realizations of the random perturbation. The program perturbation, on the other hand, produces feasible queries (up to prescribed 1\% tolerance) but at the expense of an optimality loss. As expected, the optimality loss is growing in adjacency parameter $\alpha$ as more datasets are made statistically indistinguishable in query answers, revealing the underlying cost-privacy trade-offs. Notably, in case of 14\_ieee dataset for $\alpha=10$, the chance-constrained program \eqref{OPF:sto} failed to return the solution (program is infeasible), meaning that the given level of privacy cannot be provided at the desired feasibility level.

\begin{figure}
    \FIGURE{
    \begin{subfigure}[t]{0.45\textwidth}
        \centering
        \includegraphics[width=1\textwidth]{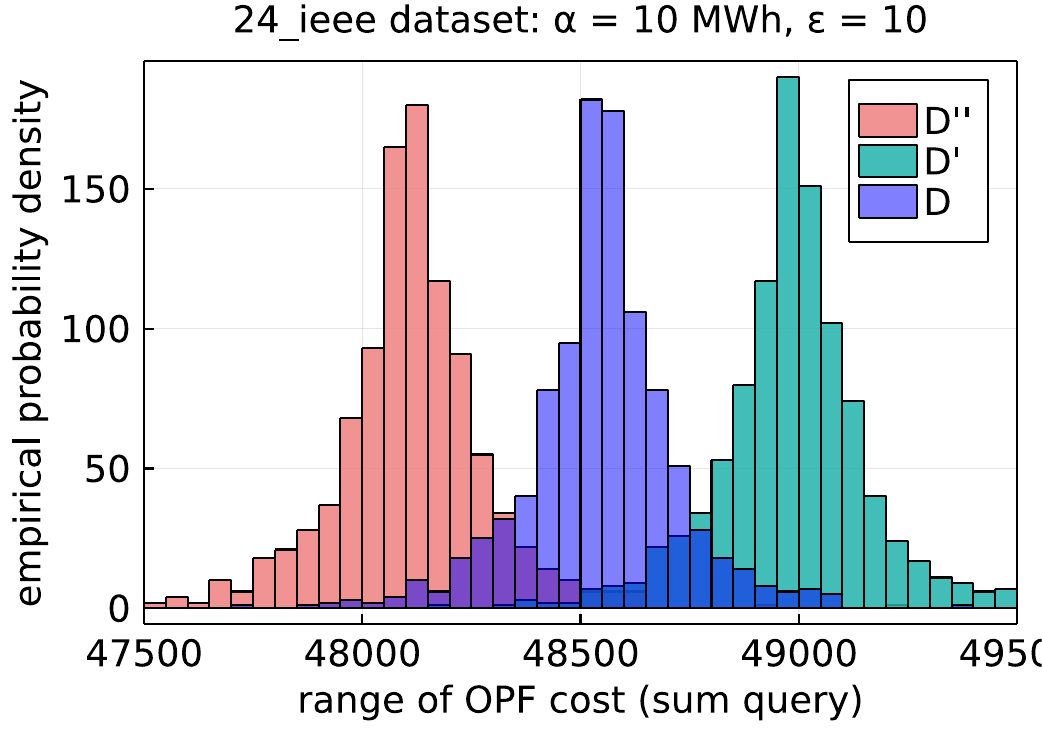}
    \end{subfigure}%
    ~ 
    \begin{subfigure}[t]{0.45\textwidth}
        \centering
        \includegraphics[width=1\textwidth]{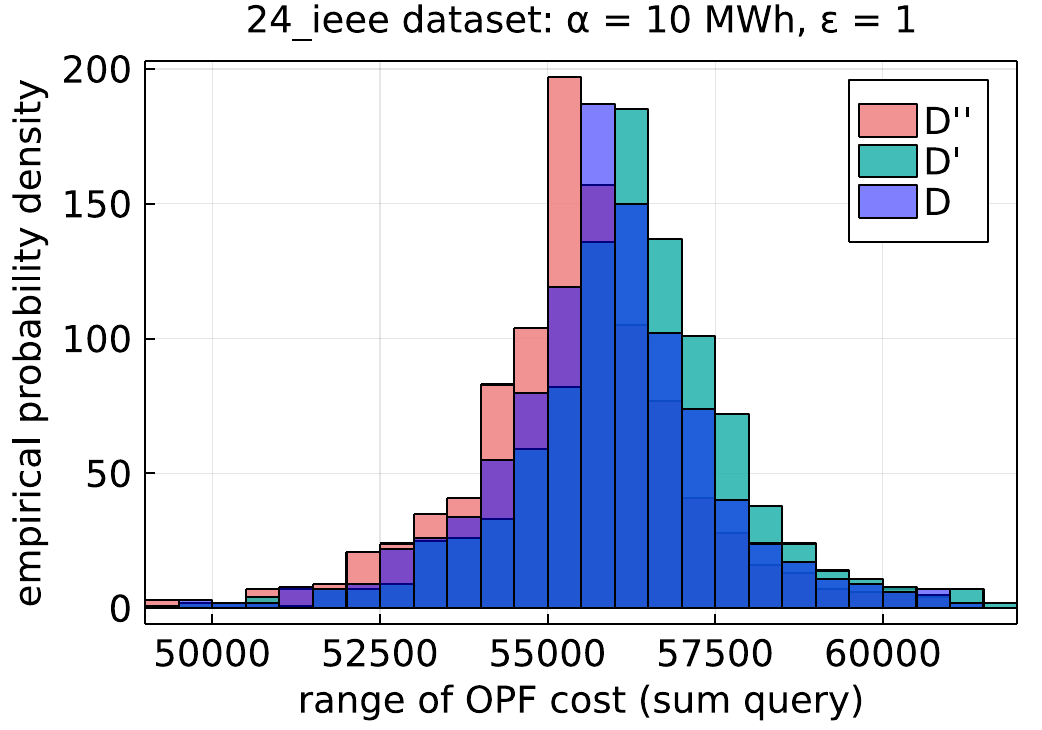}
    \end{subfigure}%
    }
    {Visualization of the probabilistic privacy guarantees under the program perturbation strategy using empirical probability density of the objective function value on three adjacent datasets.\label{fig:OPF_query_guarantee}}
    {}
\end{figure}

For the program perturbation strategy, in Figure \ref{fig:OPF_query_guarantee}, we visualize perturbed sum query distributions induced on three $\alpha-$adjacent datasets $\mathcal{D}''$, $\mathcal{D}'$ and $\mathcal{D}$, where the maximal element in demand vector $d$ is subtracted by $\alpha$, increased by $\alpha$ or left intact, respectively. For privacy parameter $\varepsilon=10$ (upper plot), the three distributions are distinguished, clearly revealing the differences in datasets. For  $\varepsilon=1$ (lower plot), all distributions are co-aligned in the query's range, thereby providing a stronger privacy protection. However, the resulting distributions are centered at a larger cost value.  

Last, we study how to control the distribution of optimality loss w.r.t. the non-private solution using CVaR model \eqref{prog:cvar}\textcolor{maincolor}{, with $S=1,000$}. Consider a query over a subset of vector $x$, such as an aggregated power generation by a particular technology, e.g., aggregated wind power generation. Table \ref{tab:CVaR} collects the mean and $\text{CVaR}_{5\%}$ values of the optimality loss induced by such a query over 50\% of generation units (chosen randomly) in 24\_ieee dataset for the fixed privacy parameter $\varepsilon=1$ and varying adjacency and CVaR parameters $\alpha$ and $q$, respectively. As privacy requirement $\alpha$ increases, the optimality loss also increases. For $q=99\%$, the average optimality loss over 5\% of the worst-case scenarios is almost twice as much as the mean value. By reducing quantile $q$, however, the worst-case suboptimality can be reduced at the expense of an increasing mean value. Eventually, the mean and $\text{CVaR}_{5\%}$ values can be matched, meaning that the same level of privacy (randomization) can be achieved at a fixed optimality loss. Hence, the curator of the private OPF problem is given a choice between playing a lottery or paying a fixed cost of privacy.

\begin{table}
\TABLE{Mean and $\text{CVaR}_{5\%}$ values of the random optimality loss of private aggregated generation query \label{tab:CVaR}}
{
\footnotesize
\begin{tabular}{ccccccccccccc}
\toprule
\multicolumn{2}{c}{\multirow{2}{*}{dataset adjacency}} & \multicolumn{11}{c}{\% of the worst-case suboptimality scenarios ($q$ in model \eqref{prog:cvar})} \\
\cmidrule(lr){3-13}
\multicolumn{2}{c}{} & 99 & 90 & 80 & 70 & 60 & 50 & 40 & 30 & 20 & 10 & 1 \\
\midrule
\multirow{2}{*}{$\alpha=15$ MWh} & mean & 1.56 & 1.56 & 1.56 & 1.70 & 1.70 & 1.70 & 1.72 & 1.72 & 1.72 & 1.72 & 1.72 \\
 & $\text{CVaR}_{5\%}$ & 2.36 & 2.36 & 2.36 & 1.78 & 1.78 & 1.78 & 1.72 & 1.72 & 1.72 & 1.72 & 1.72 \\
\midrule
\multirow{2}{*}{$\alpha=25$ MWh} & mean & 4.90 & 4.90 & 4.90 & 5.39 & 5.39 & 5.54 & 5.54 & 5.54 & 5.54 & 5.54 & 5.54 \\
 & $\text{CVaR}_{5\%}$ & 8.06 & 8.06 & 8.06 & 5.96 & 5.96 & 5.54 & 5.54 & 5.54 & 5.54 & 5.54 & 5.54 \\
\midrule
\multirow{2}{*}{$\alpha=35$ MWh} & mean & 8.54 & 8.54 & 8.81 & 9.03 & 9.03 & 9.03 & 9.09 & 9.44 & 9.44 & 9.44 & 9.86 \\
 & $\text{CVaR}_{5\%}$ & 14.28 & 14.28 & 11.83 & 10.89 & 10.89 & 10.89 & 10.77 & 10.33 & 10.33 & 10.33 & 10.31 \\
\bottomrule
\end{tabular}
}
{}
\end{table}

\subsection{Private Quadratic Programming: Support Vector Machines (SVM)}\label{subsec:SVM}
Consider a dataset $(x_{1},y_{1}),\dots,(x_{m},y_{m})$ of $m$ data points, where $x_{i}\in\mathbb{R}^{n}$ is a vector of features and $y_{i}\in\{-1,1\}$ is an associated binary label (class). The SVM goal is to compute the maximum-margin hyperplane that separates data points from different classes so as to maximize the distance between the hyperplane and the nearest point. Formally, SVM optimizes hyperplane $w^{\top}x_{i}-b$ with parameters $w\in\mathbb{R}^{n}$ and $b\in\mathbb{R}$, so that any data point $\tilde{x}$ is classified according to the rule $$\text{sign}[w^{\star\top}\tilde{x}-b^{\star}],$$ i.e., depending on which side of the hyperplane the point falls. The optimal hyperplane parameters are obtained solving the following convex optimization program \citep{cortes1995support}:
\begin{subequations}\label{SVM:det}
\begin{align}
    \minimize{b,w,z}\quad&\lambda\norm{w}^{2} + \textstyle\frac{1}{m}\mathbb{1}^{\top}z\\
    \st\quad&y_{i}(w^{\top}x_{i}-b)\geqslant1-z_{i},\quad z_{i}\geqslant0,\quad\forall i=1,\dots,m,
\end{align}
\end{subequations}
with regularization parameter $\lambda$ and slack variable $z\in\mathbb{R}^{m}$ to ensure feasibility when the data points are not fully separable by the hyperplane. Notably, the normal vector $w^{\star}$ of the hyperplane is unique on a particular dataset, and releasing the classification rule leads to privacy breaches. 

To release the hyperplane parameters privately, we introduce perturbation $\bs{\zeta}\in\mathbb{R}^{n+1}$ ($n$ elements of $w$ and 1 of $b$) and put forth the following stochastic counterpart of program \eqref{SVM:det}:
\begin{subequations}\label{SVM:pertrubed}
\begin{align}
    \minimize{\overline{b},\overline{w},\overline{z},B,W,Z}\quad&\mathbb{E}\left[\lambda\norm{\overline{w}+W\bs{\zeta}}^{2} + \textstyle\frac{1}{m}\mathbb{1}^{\top}(\overline{z}+Z\bs{\zeta})\right]\label{SVM_pertrubed_obj}\\
    \st\quad&\text{Pr}\left[
    \begin{array}{l}
    y_{i}((\overline{w}+W\bs{\zeta})^{\top}x_{i}-(\overline{b}+B\bs{\zeta}))\geqslant1-\overline{z}_{i}-Z_{i}^{\top}\bs{\zeta},\\
    \overline{z}_{i}-Z_{i}^{\top}\bs{\zeta}\geqslant0,\quad\forall i=1,\dots,m,
    \end{array}
    \right]\geqslant 1-\eta,
    \quad 
    \begin{bmatrix}
    W \\
    B
    \end{bmatrix} = I,
\end{align}
\end{subequations}
which minimizes the expected objective function value subject to a chance constraint on perturbed SVM constraints, where linear decision rules are enforced for each decision variable. The recourse decisions $B$ and $W$ are constrained to respect the identity query constraints, as in Example \ref{example:identity_query}. Once problem \eqref{SVM:pertrubed} is solved to optimality, the differentially private hyperplane release includes the intercept $\overline{b}^{\star} + B^{\star}\hat{\zeta}$ and the normal vector $\overline{w}^{\star} + W^{\star}\hat{\zeta}$, where $\hat{\zeta}$ is a random realization of $\bs{\zeta}$.

\textcolor{maincolor}{Towards a tractable reformulation, the linear term in \eqref{SVM_pertrubed_obj} reformulates as $\mathbb{E}\left[\frac{1}{m}\mathbb{1}^{\top}(\overline{z}+Z\bs{\zeta})\right] = \frac{1}{m}\mathbb{1}^{\top}\overline{z}$ due to zero mean of perturbation $\bs{\zeta}$, and the  quadratic term as $\mathbb{E}\left[\lambda\norm{\overline{w}+W\bs{\zeta}}^{2}\right] = \norm{\overline{w}}^2 + \text{Tr}\left[W\Sigma W^{\top}\right],$ where $\Sigma$ is a covariance matrix.
The joint chance constraint is split into $2m$ individual chance constraints with constraint violations set uniformly as $\overline{\eta}=\eta/(2m)$. They are then reformulated as in Remark \ref{remark:lin_con_cc}, allowing for the analytic SOCP reformulation of the entire program.}

\subsubsection{SVM on Synthetic Data.}\label{subsubsec:SVM_example}

\begin{figure}
\FIGURE{
\includegraphics[width=0.95\textwidth]{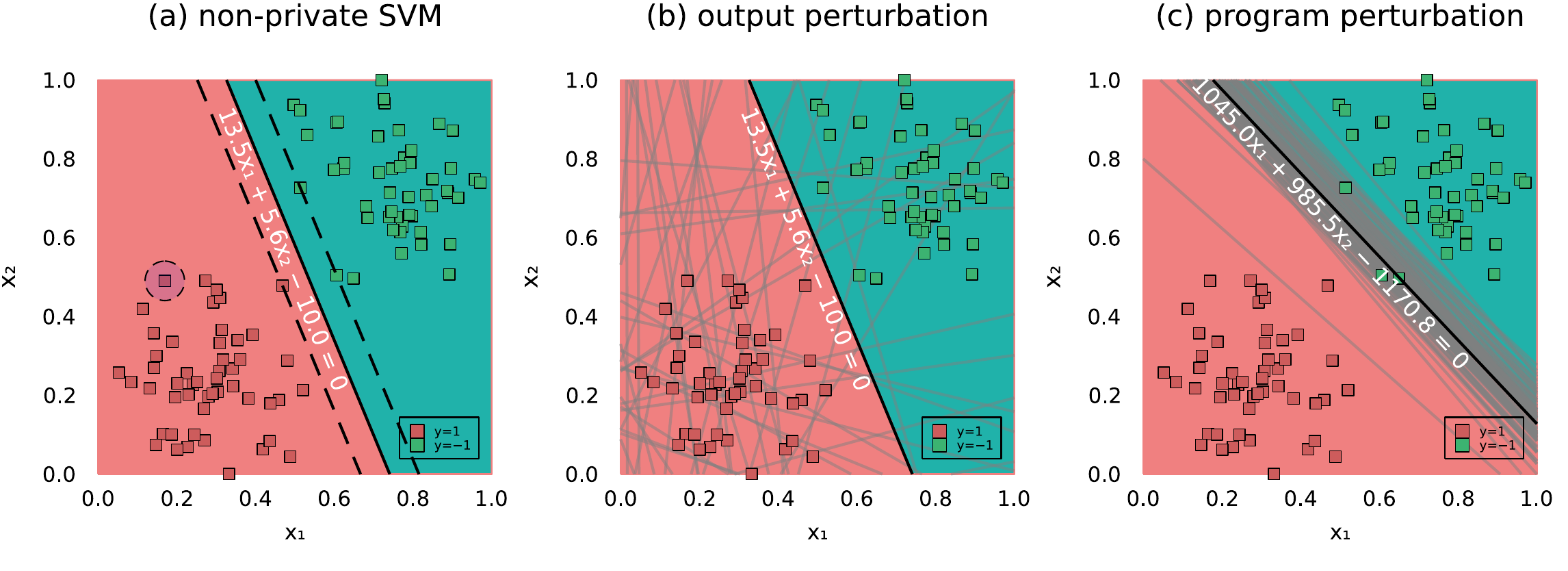}}
{
Privacy-preserving SVM classification on synthetic data: (a) Problem data and deterministic, non-private SVM solution, (b) Output perturbation with the privacy-preserving Laplacian noise, and (c) Proposed program perturbation with the same noise as in the output perturbation. The solid black line separates data points and the dashed lines on the left figure illustrate the data separating margin. The gray solid lines depict random realizations of perturbed separating hyperplanes. The input perturbation strategy has not been implemented here to keep the original data points unchanged.\label{fig:SVM_example}
}{}
\end{figure}

We generate a classification dataset with 100 training data points from distribution
$\mathcal{N}(
\begin{bsmallmatrix}
1\\1
\end{bsmallmatrix},
\begin{bsmallmatrix*}[l]
0.5 & 0.0\\ 0.0 & 0.5
\end{bsmallmatrix*})$ for class $y=1$ and distribution
$\mathcal{N}(
\begin{bsmallmatrix}
3\\3
\end{bsmallmatrix},
\begin{bsmallmatrix*}[l]
0.5 & 0.0\\ 0.0 & 0.5
\end{bsmallmatrix*})$
for class $y=-1$ in equal proportion, and then normalize them using a min-max normalization. After solving program \eqref{SVM:det} with small regularization parameter $\lambda=10^{-5}$, the SVM hyperplane separates the training dataset as shown in Figure \ref{fig:SVM_example}(a). The classification accuracy on $10^3$ testing data points amounts to 99.4\%.

The privacy goal here is to make the training dataset indistinguishable from any adjacent dataset when releasing hyperplane parameters. We assume such a dataset universe $\mathbb{D}$, where $i^{\text{th}}$ training data point realizes according to the law $(r\sin(t) + x_{1i}, r\cos(t) + x_{2i}),$ with $r\sim U(0,0.05)$ and $t\sim U(0,2\pi)$, i.e., from the scaled circle area centered at $(x_{1i},x_{2i})$, as shown in Figure \ref{fig:SVM_example}(a) for a selected data point, and we let all training datasets from this universe to be adjacent. The sensitivity $\Delta_{p}$ of hyperplane parameters to datasets is evaluated at 21.8 using Algorithm \ref{alg:sens} with adjacency parameter $\alpha=\infty$, norm order $p=1$, and sample size requirement $S=99$ (for $\gamma=\beta=0.1$ of Proposition \ref{prop:sample_size}).

To provide a pure $1-$differential privacy ($\varepsilon=1$) to training datasets, we use Laplace perturbation $\bs{\zeta}\sim\text{Lap}(1/21.8)$ and first deploy the output perturbation strategy (Definition \ref{def:OP}), i.e., by perturbing the solution of program \eqref{SVM:det}. Figure \ref{fig:SVM_example}(b) displays 100 random hyperplanes (in gray) resulting from the output perturbation. The data separations are of poor quality with a classification accuracy of 51.2\% on average with a standard deviation of 11.6\% on the testing dataset. 

Next, we deploy the proposed program perturbation strategy by perturbing the nominal solution of chance-constrained program \eqref{SVM:pertrubed} with the same Laplace noise as in the output perturbation while requiring solution feasibility with a probability of at least $1-\eta=95\%$. Figure \ref{fig:SVM_example}(c) displays 100 random hyperplanes resulting from the program perturbation. Observe, that the nominal hyperplane parameters exceed those of the deterministic solution by 2 orders of magnitude (e.g., 1045.0 versus 13.5 for parameter $w_{1}$). Thus, the hyperplane becomes less sensitive to perturbation, which results in a more stable distribution of the random hyperplanes around the nominal solution. The strategy demonstrates a classification accuracy of 97.6\% on average with a standard deviation of 1.7\% on the testing dataset. Hence, the desired privacy protection of the training dataset is achieved while maintaining the classification accuracy close to that of the non-private solution.

\subsubsection{SVM for OPF Feasibility Classification.}
For a given vector of demands $d_{i}$, the instance of the OPF problem in \eqref{OPF:det} can be classified as feasible $(y_{i}=1)$ or infeasible $(y_{i}=-1)$, depending on the solution status of problem \eqref{OPF:det}. We use differentially private SVM to classify OPF feasibility without disclosing the labeled demand dataset. For each selected network from \citep{coffrin2018powermodels}, we provide three training datasets in Table \ref{tab:OPF_SVM_results} depending on how labels are balanced (see the $2^{\text{nd}}$ column). All datasets include $m=4500$ rows, each including $n$ features, where $n$ amounts to non-zero demands in vector $d$ (see the $3^{\text{rd}}$ column). 

\begin{table}
\TABLE
{Summary of $(\varepsilon,\delta)-$differentially private OPF feasibility classification. The mean classification accuracy is supplemented with the standard deviation in parentheses. Here, $\overline{\eta}_{i}$ indicates the tolerance to each SVM constraint violation (formulated for each data point $i=1,\dots,m$). 
\label{tab:OPF_SVM_results}
}
{
\footnotesize
\begin{tabular}{lccccccc}
\toprule
\multicolumn{1}{c}{\multirow{3}{*}{dataset}} & \multirow{2}{*}{\begin{tabular}[c]{@{}c@{}}dataset \\ balance (\%)\end{tabular}} & \multirow{3}{*}{\begin{tabular}[c]{@{}c@{}} \# of \\ features \\ (demands)\end{tabular}} & \multicolumn{5}{c}{mean classification accuracy (\%)} \\
\cmidrule(lr){4-8}
\multicolumn{1}{c}{} &  &  & \multirow{2}{*}{\begin{tabular}[c]{@{}c@{}}non-private\\ solution\end{tabular}} & \multirow{2}{*}{\begin{tabular}[c]{@{}c@{}}output \\ perturbation\end{tabular}} & \multicolumn{3}{c}{program perturbation} \\
\cmidrule(lr){2-2}\cmidrule(lr){6-8}
\multicolumn{1}{c}{} & $-1$ / $+1$ &  &  &  & $\overline{\eta}_{i}=5\%$ & $\overline{\eta}_{i}=1\%$ & $\overline{\eta}_{i}=0.01\%$ \\
\midrule
5\_pjm\_1 & 4.3\textcolor{white}{0} / 95.7 & $|$ & 99.6 & 79.7 (22.5) & 95.6 (2.0) & 95.0 (1.3) & 94.3 (0.4) \\
5\_pjm\_2 & 26.2 / 73.8 & 3 & 96.2 & 75.8 (12.8) & 93.4 (3.0) & 94.3 (2.2) & 95.3 (1.3) \\
5\_pjm\_3 & 50.6 / 49.4 & $|$ & 95.6 & 79.1 (12.2) & 93.6 (2.4) & 94.3 (1.6) & 95.1 (0.8) \\
\midrule
14\_ieee\_1 & 5.4\textcolor{white}{0} / 94.6 & $|$ & 98.6 & 76.6 (26.3) & 93.9 (2.1) & 93.6 (1.1) & 93.2 (0.2) \\
14\_ieee\_2 & 26.3 / 73.7 & 11 & 98.4 & 71.5 (18.3) & 93.6 (3.5) & 94.4 (2.4) & 94.9 (1.1) \\
14\_ieee\_3 & 47.7 / 52.3 & $|$ & 98.6 & 69.7 (14.5) & 91.4 (3.6) & 92.1 (2.4) & 91.8 (1.1) \\
\midrule
24\_ieee\_1 & 5.1\textcolor{white}{0} / 94.9 & $|$ & 96.6 & 77.4 (24.9) & 94.7 (2.7) & 94.8 (1.2) & 94.6 (0.1) \\
24\_ieee\_2 & 22.9 / 77.1 & 17 & 92.2 & 72.9 (16.9) & 89.3 (3.9) & 90.4 (2.9) & 91.8 (1.7) \\
24\_ieee\_3 & 48.5 / 51.5 & $|$ & 89.4 & 66.9 (10.4) & 85.1 (4.4) & 86.3 (3.5) & 87.3 (2.1) \\
\midrule
57\_ieee\_1 & 4.1\textcolor{white}{0} / 95.9 & $|$ & 99.0 & 69.5 (34.5) & 95.4 (6.4) & 96.2 (2.0) & 96.4 (0.0) \\
57\_ieee\_2 & 23.8 / 76.2 & 42 & 96.0 & 64.0 (25.3) & 91.3 (3.3) & 92.0 (2.4) & 92.4 (1.6) \\
57\_ieee\_3 & 47.1 / 52.9 & $|$ & 97.4 & 62.1 (14.7) & 90.2 (4.0) & 91.0 (2.6) & 91.3 (1.2) \\
\midrule
89\_pegase\_1 & 3.1\textcolor{white}{0} / 96.9 & $|$ & 97.4 & 51.1 (44.7) & 93.1 (15.6) & 95.8 (5.4) & 96.2 (0.0) \\
89\_pegase\_2 & 25.4 / 74.6 & 35 & 86.8 & 50.4 (25.0) & 75.5\textcolor{white}{0} (6.4) & 75.7 (2.0) & 75.6 (0.0) \\
89\_pegase\_3 & 50.5 / 49.5 & $|$ & 82.2 & 50.3\textcolor{white}{0} (3.4) & 73.8\textcolor{white}{0} (7.5) & 75.6 (5.9) & 77.1 (4.2)\\
\bottomrule
\end{tabular}
}
{}
\end{table}

We provide $(\varepsilon,\delta)-$differential privacy (Theorem \ref{th:eps_delta_dp_id}), with $\varepsilon=1$ and $\delta=\frac{1}{m}=2.2\times10^{-4}$, for $\alpha=\infty$ classification datasets from the universe, where each data point (demand) varies in the range $\pm1\%$ of its nominal value. The rest of the settings are the same as in Section \ref{subsubsec:SVM_example}. The output perturbation results are summarized in the $4^{\text{th}}$ column and demonstrate a substantial reduction in classification accuracy across all networks. The program perturbation results are provided in the remaining columns for varying tolerance to SVM constraint violations. As program \eqref{SVM:pertrubed} minimizes the expected value of the perturbed objective function (the proxy of classification accuracy) and internalizes feasibility criteria, the program perturbation demonstrates a substantially smaller loss in the classification accuracy. Moreover, tightening the constraint violation tolerance reduces the standard deviation of the classification accuracy. This is because chance constraints ``shrink" the feasible region for the nominal solution $\overline{w}$ and $\overline{b}$, making it more robust to perturbation.

\subsection{Private Quadratic Programming: Monotonic Regression}\label{subsec:Regression}

Consider the problem of fitting a smooth function $h(x)$ to a set of observations $(x_{1},y_{1}),\dots,(x_{n},y_{n})$, subject to the constraint that $h$ must be monotone in the range of $x$. Such constraints originate in regression problems, where the underlying phenomenon is known to be monotone, e.g., wind power curve fitting \citep{mehrjoo2020wind}. The function is modeled as $h(x)=w^{\top}\varphi(x)$, where $w\in\mathbb{R}^{m}$ is a vector of weights and $\varphi(x)\in\mathbb{R}^{m}$ is a high-dimensional basis function. Function $h(x)$ fits data when the regularized regression loss $\sum_{i=1}^{n}\norm{y_{i} - w^{\top}\varphi(x_{i})} + \lambda\norm{w}$ is minimized subject to constraint $Cw\geqslant\mathbb{0}$, where matrix $C\in\mathbb{R}^{p\times m}$ encodes monotonic conditions at $p$ selected points, i.e., its $i^{\text{th}}$ row contains partial derivatives $\varphi'(u_{i})$ evaluated at selected point $u_{i}$ in the range of $x$. 

The privacy goal is to make regression weights $w$ statistically similar on adjacent regression datasets. The application of the standard output perturbation strategy $w^{\star} + \hat{\zeta}$, however, does not guarantee that the perturbed weights retain monotonic conditions, i.e., that $C(w^{\star} + \hat{\zeta})\geqslant\mathbb{0}$. To guarantee feasibility, we redefine the weights as linear decision rule $w(\bs{\zeta}) = \overline{w} + W\bs{\zeta}$, where $\overline{w}\in\mathbb{R}^{m}$ and $W\in\mathbb{R}^{m\times m}$ are optimization variables, and then solve the following stochastic program
\begin{subequations}\label{MCF:pertrubed}
\begin{align}
    \minimize{\overline{w},W}\quad&\mathbb{E}\left[\sum_{i=1}^{n}\norm{y_{i} - (\overline{w} + W\bs{\zeta})^{\top}\varphi(x_{i})} + \lambda\norm{\overline{w} + W\bs{\zeta}}\right]\\
    \st\quad&\text{Pr}\left[
    C(\overline{w} + W\bs{\zeta})\geqslant\mathbb{0}
    \right]\geqslant1-\eta,\quad W = I,\label{MCF:cc}
\end{align}
\end{subequations}
which minimizes the expected regression loss, while requiring constraint satisfaction with probability $(1-\eta)$ and enforcing the identity query constraints to guarantee noise independence, as in Example \ref{example:identity_query}. Unlike the deterministic problem, this stochastic program optimizes the mean weights $\overline{w}$ anticipating the impact of the perturbation on the probability of monotonic constraint satisfaction. Once program \eqref{MCF:pertrubed} is solved to optimality, the differentially private regression weights release is $\overline{w}^{\star} + \hat{\zeta}$, where $\hat{\zeta}$ is a random realization of $\bs{\zeta}$. \textcolor{maincolor}{The tractable analytic reformulation of problem \eqref{MCF:pertrubed} is achieved in the same manner as that of problem \eqref{SVM:pertrubed} in Section \ref{subsec:SVM}.}

\subsubsection{Monotonic Regression on Synthetic Data.}\label{subsubsec:monoton}

\begin{figure}
    \FIGURE{
    \begin{subfigure}[t]{0.33\textwidth}
        \centering
        \includegraphics[width=1\textwidth]{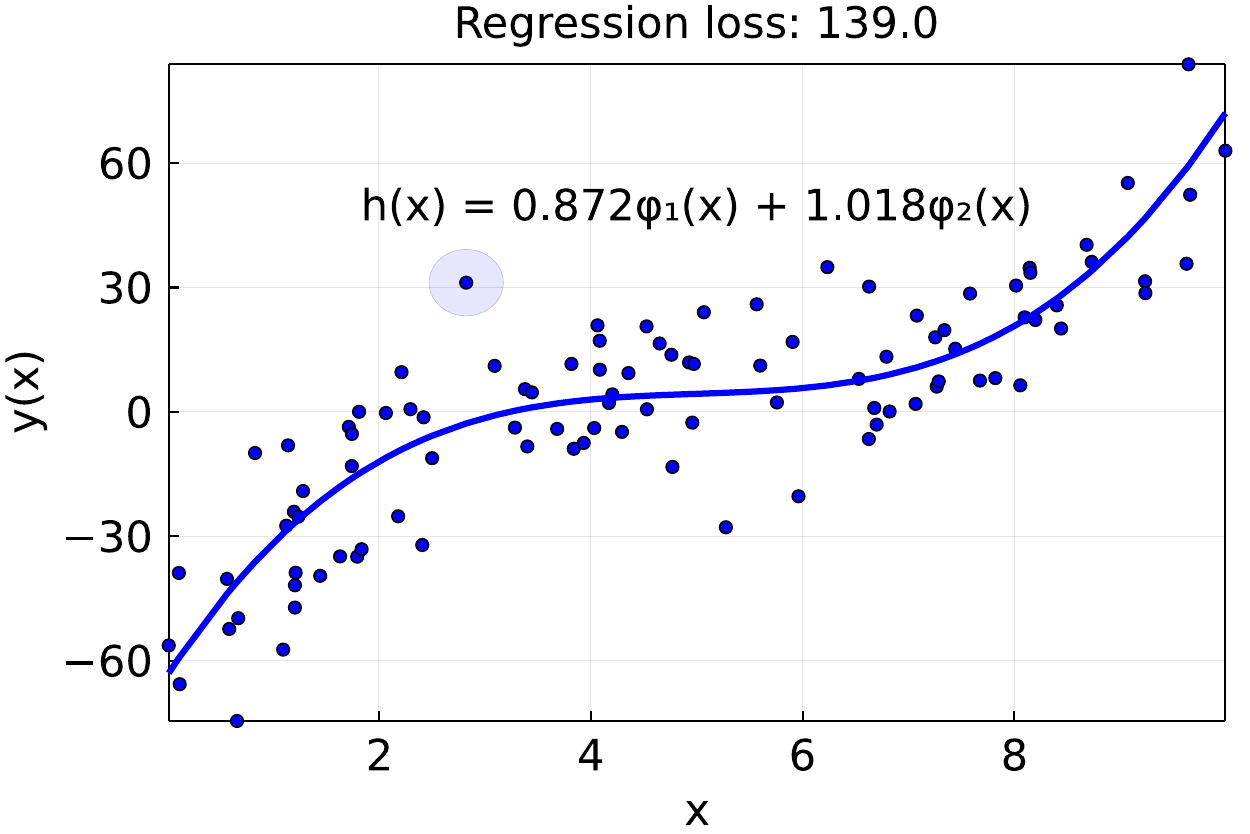}
        \caption{Baseline, non-private model}\label{fig:reg_base}
    \end{subfigure}%
    ~ 
    \begin{subfigure}[t]{0.33\textwidth}
        \centering
        \includegraphics[width=1\textwidth]{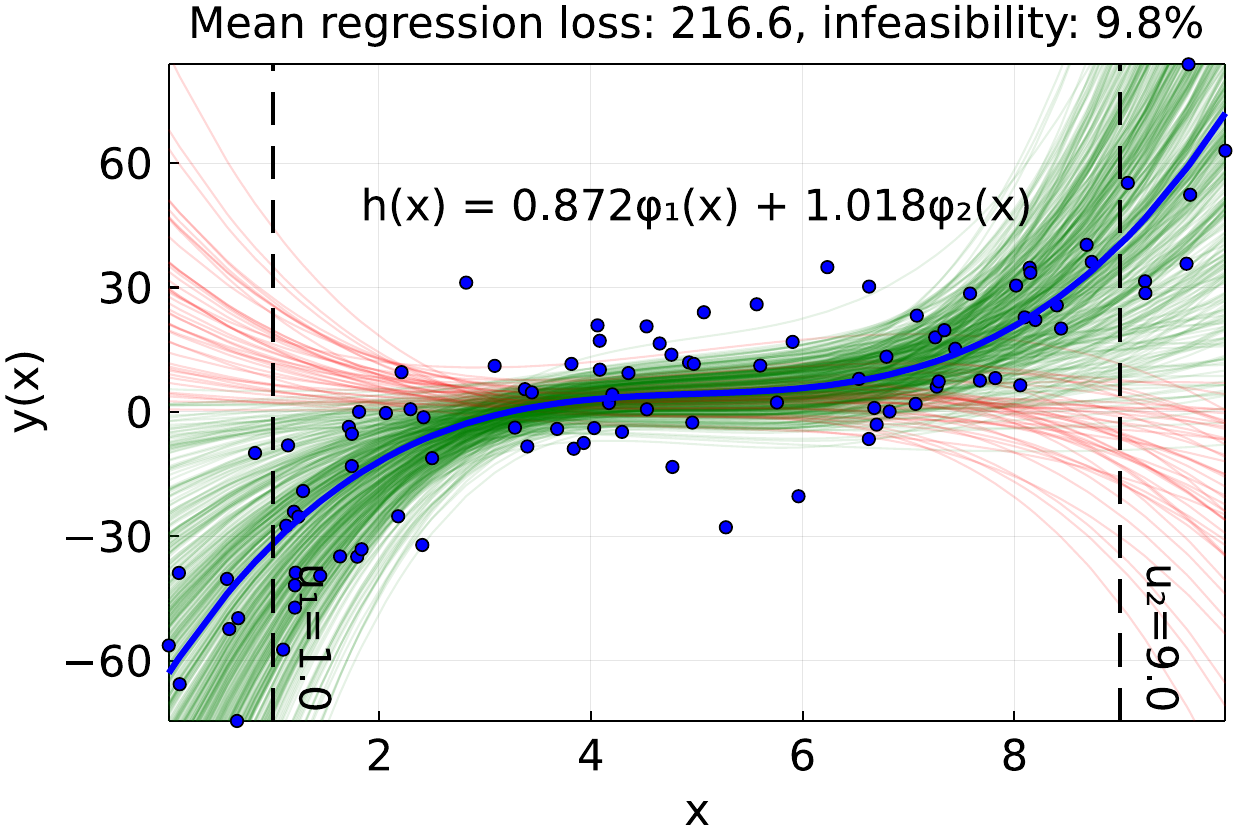}
        \caption{Output perturbation}\label{fig:reg_op}
    \end{subfigure}%
    \begin{subfigure}[t]{0.33\textwidth}
        \centering
        \includegraphics[width=1\textwidth]{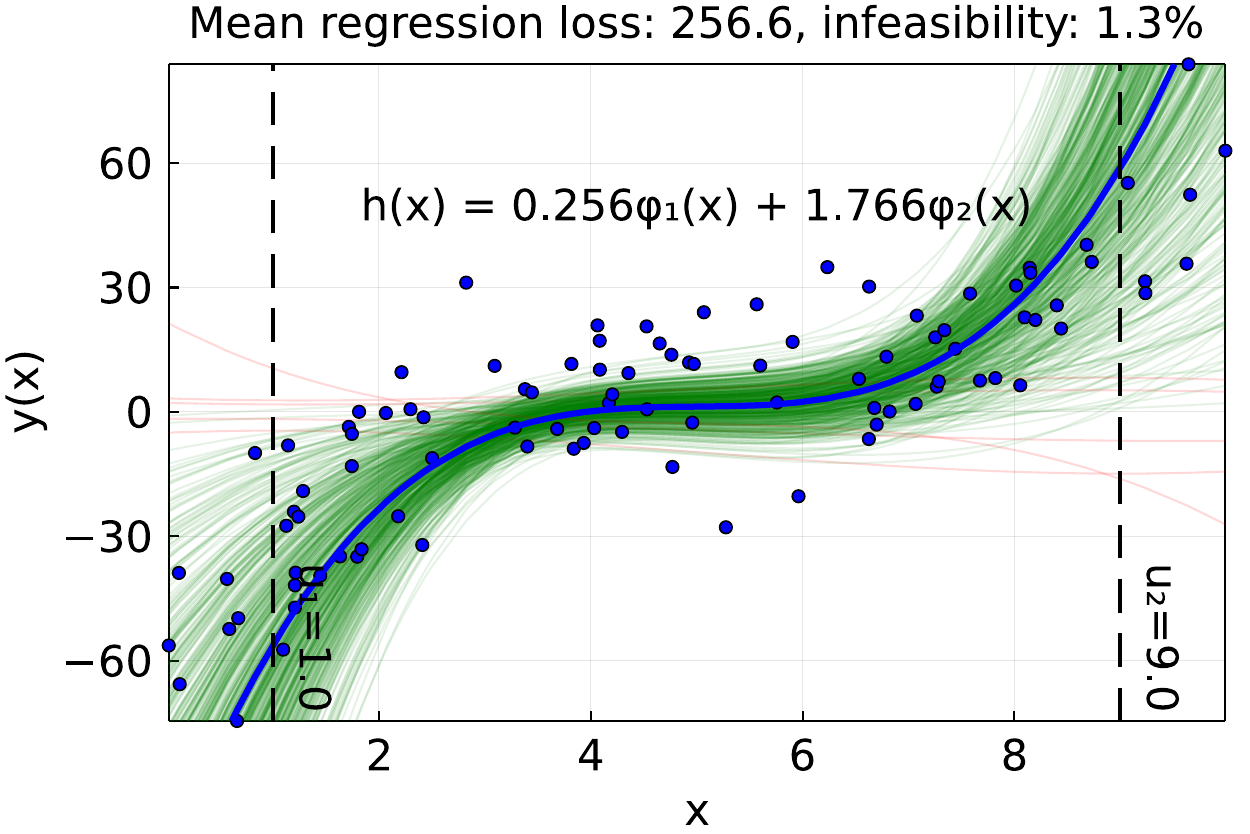}
        \caption{Program perturbation}\label{fig:reg_pp}
    \end{subfigure}%
    }
    {Differentially private monotonic regression: (a) Problem data, the baseline regression model and dataset universe $\mathbb{D}$ depicted as a circle domain around a single selected point; (b) Output perturbation results: 500 perturbed regression models depicted in green and red depending on whether they satisfy the monotone increasing conditions at points $u_{1}=1$ and $u_{2}=9$; (c) Program perturbation results: the mean regression model deviates from the baseline, but the perturbation of its weights does not violate monotone increasing conditions with a very high probability.}
    {}
\end{figure}

We illustrate the private monotonic regression by approximating the following function:
\begin{align*}
y(x) = \mathbb{1}^{\top}
\underbrace{\begin{bmatrix}
x\\
\frac{1}{2}(x-5)^3
\end{bmatrix}}_{\varphi(x)}
+z,\quad z\sim N(0,15),
\end{align*}
with the linear and cubic basis functions and some random noise $z$. The regression dataset is generated from this function using 100 independent variables $x$ drawn from the uniform distribution $U(0,10)$. Figure \ref{fig:reg_base} illustrates the dataset and the baseline regression model.

We assume a dataset universes $\mathbb{D}$, where the $i^{\text{th}}$ data point realizes according to the law $(0.35r\cos(t)+x_{i},8r\sin(t)+y_{i}),$ where $r\sim U(0,1)$ and $t\sim U(0,2\pi)$, i.e., from the scaled circle domains centered at $(x_{i},y_{i})$, as further shown in Figure \ref{fig:reg_base} for a selected data point. In this example, we provide $(1,0.01)-$differential privacy for any dataset from this universe, i.e., $\varepsilon=1,\delta=0.01$ and $\alpha=\infty$. The sensitivity $\Delta_{2}$ of the baseline model (regression weights) to datasets evaluates at 0.46 using Algorithm \ref{alg:sens} with the sample size requirement $S=199$ (for $\gamma=0.5,\beta=0.1$ of Proposition \ref{prop:sample_size}). 

We first deploy the output perturbation strategy using the Gaussian mechanism, i.e., by perturbing the weights of the baseline model with the Gaussian noise $N(0,\sqrt{2\text{ln}(1.25/0.01)}0.46)$ as per Theorem \ref{th:eps_delta_dp_id}, and illustrate the results in  Figure \ref{fig:reg_op}. Observe, that the output perturbation is likely to violate the original monotonic conditions of the baseline model. To measure violations, we select two points $u_{1}=1$ and $u_{2}=9$ in the range of $x$ and build the constraint matrix
\begin{align*}
C=
\begin{bmatrix}
\partial\varphi_{1}(u_{1})/\partial u_{1} & \partial\varphi_{2}(u_{1})/\partial u_{1}\\
\partial\varphi_{1}(u_{2})/\partial u_{2}& \partial\varphi_{2}(u_{2})/\partial u_{2}
\end{bmatrix}
=
\begin{bmatrix}
1.00 & 36.75\\
1.00 & 48.00
\end{bmatrix},
\end{align*}
and we record a violation whenever $C(w^{\star}+\hat{\zeta})\geqslant0$ does not hold under perturbation. For the output perturbation, this constraint is violated in 9.8\% instances across 500 perturbed models shown in either green (monotone increasing models) or red (monotone decreasing models) in Figure \ref{fig:reg_op}. 

Next, we deploy the proposed program perturbation strategy by first solving chance-constrained program \eqref{MCF:pertrubed} and then perturbing its mean solution $\overline{w}$. Here, we require the monotonic constraint satisfaction with a probability of at least $97\%$ ($\eta$ is set to 3\%). The corresponding results are depicted in Figure \ref{fig:reg_pp}, where the mean regression model deviates from the baseline model, but the constraint violation probability on 500 samples is as small as 1.3\%, i.e., below the prescribed maximum violation probability of 3\%. This feasibility guarantee, however, comes at the expense of an increasing regression loss, which exceeds that of the output perturbation by 18.4\% in expectation, thus revealing the inherent trade-offs between feasibility and the goodness of fit.  

\subsubsection{Monotonic Regression for Wind Power Curve Fitting.}

A wind turbine has a unique power curve characteristic that explains how power generation depends on weather inputs, e.g., wind speed. This curve is monotonically increasing within medium-range wind speeds, e.g., 3 to 10 m/s. Since the theoretical power curve provided by manufactures does not internalize specifics of operational and geographic conditions, it is common to estimate this curve from historical records. These records often include private information (e.g., operational patterns, anomalies such as icing) which wind turbine operators are unwilling to disclose. It is thus relevant to produce physically-meaningful (monotone) wind power curves that obfuscate this private information. 

We use  Renewable Ninja database \citep{staffell2016using} to collect normalized, theoretical wind power curves from several manufactures (see the list of selected wind turbines in Table \ref{tab:wind_power_curve_results}). To simulate the real-world historical records, every data point is perturbed with the Gaussian noise from $N(0,0.1)$ and post-processed to ensure consistency (e.g., w.r.t. minimum and maximum power bounds of each turbine). Each data point $(y_{i},x_{i})$ includes wind power record $y_{i}$ and corresponding wind speed $x_{i}$. To approximate the curve of each turbine, we use $m=4$ radial basis functions of the form $\varphi_{i}(x) = \sqrt{1+(\mu_{i}-x)^2},\forall i=1,\dots,m,$ respectively centered at $\mu=\{3,7,11,15\}$ m/s. The remaining parameters are the same as in Section \ref{subsubsec:monoton}. 

\begin{figure}
    \FIGURE{
    \begin{subfigure}[t]{0.45\textwidth}
        \centering
        \includegraphics[width=1\textwidth]{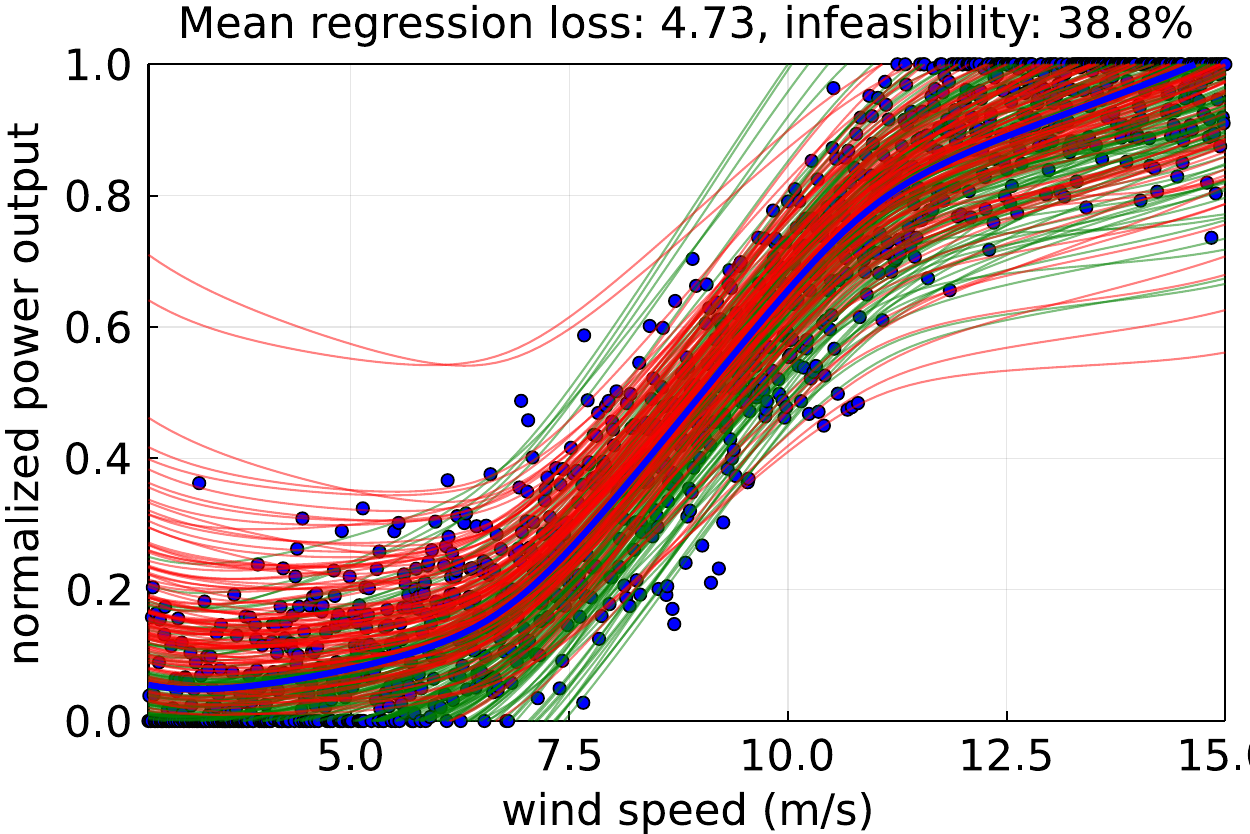}
    \end{subfigure}%
    ~ 
    \begin{subfigure}[t]{0.45\textwidth}
        \centering
        \includegraphics[width=1\textwidth]{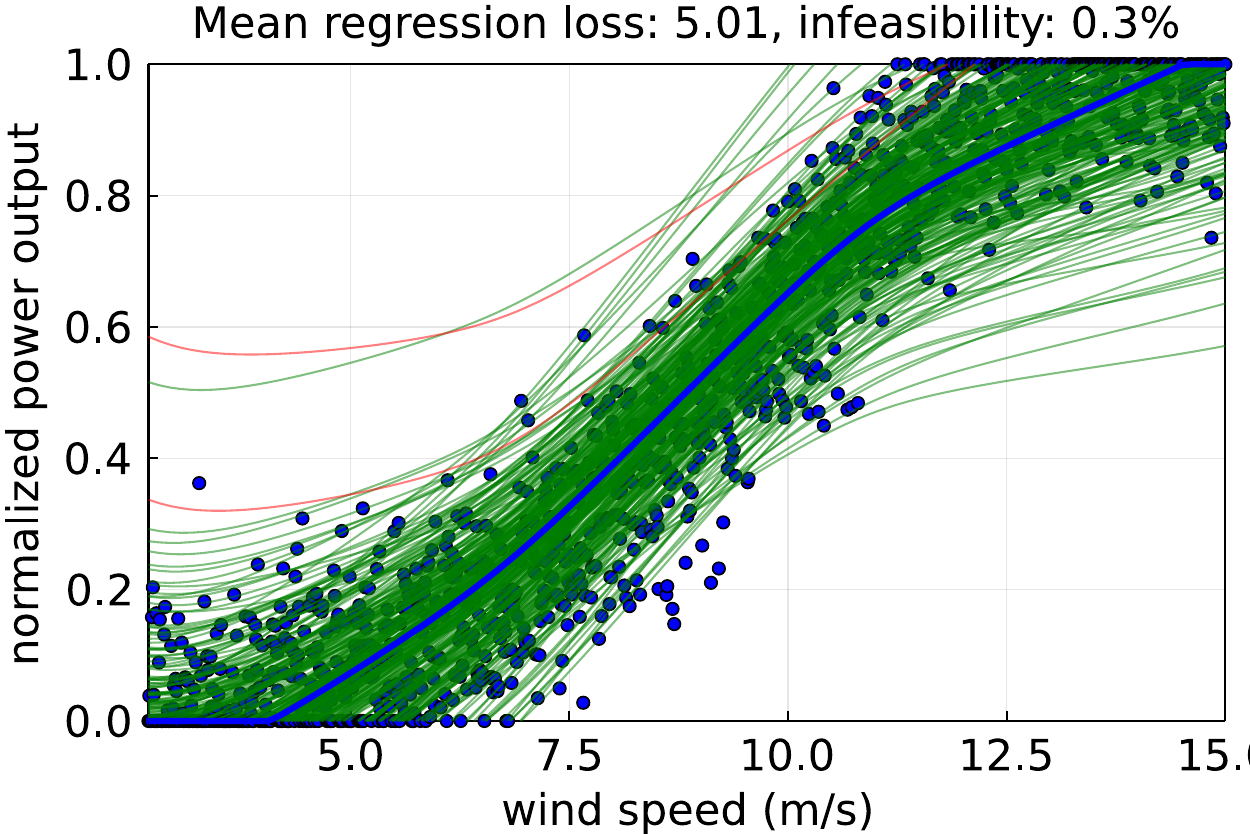}
    \end{subfigure}%
    }
    {Differentially private wind power curve realizations for Alstom.Eco.80 turbine under output (left) and program (right) perturbation strategies for dataset variation $\pm2.50\%$.\label{fig:Alstom}}
    {}
\end{figure}

Consider a dataset universe, where the power output of each data point varies within the range $\pm1\%$, $\pm1.75\%$  or $\pm2.5\%$ of the original value. As this range increases, the range of datasets that are required to be statistically indistinguishable (in the model weights query) also increases, hence improving privacy. To ensure monotonic properties, we construct matrix $C$ in \eqref{MCF:cc} by selecting $p=10$ points at random within the range from 3 to 10 m/s. The results of the output and program perturbation strategy are reported in Table \ref{tab:wind_power_curve_results}. The two strategies result in a larger regression loss than the non-private solution, and the loss increases as the domain of dataset universe increases, revealing the trade-offs between the privacy level and the model accuracy. The output perturbation, however, systematically  fails to deliver physically meaningful power curves that preserve monotonic proprieties of the non-private solution. The program perturbation, on the other hand, identifies the adjusted model wights whose perturbation is guaranteed to be feasible with a high probability. Notably, for Vestas, Siemens, Enercon and GE turbines, such an adjustment is marginal (note the marginal difference in the expected regression loss), but the gain in model feasibility is significant. The output and program perturbation strategies are visualized in Figure \ref{fig:Alstom} for a selected turbine. 

\begin{table}
\TABLE
{Regression loss and monotonic constraints violation (in percentage) of differentially private wind power curve fitting. \label{tab:wind_power_curve_results}}
{
\footnotesize
\begin{tabular}{llccccccc}
\toprule
\multicolumn{1}{l}{\multirow{3}{*}{turbine}} & \multirow{3}{*}{\begin{tabular}[l]{@{}l@{}}perturb. \\ strategy\end{tabular}} & \multicolumn{1}{c}{\multirow{3}{*}{\begin{tabular}[l]{@{}l@{}}non-private\\ reg. loss\end{tabular}}} & \multicolumn{6}{c}{datasets variation (\%)} \\
\multicolumn{1}{c}{} &  & \multicolumn{1}{c}{} & \multicolumn{2}{c}{$\pm1.00\%$} & \multicolumn{2}{c}{$\pm1.75\%$} & \multicolumn{2}{c}{$\pm2.50\%$} \\
\cmidrule(lr){4-5}\cmidrule(lr){6-7}\cmidrule(lr){8-9}
\multicolumn{1}{c}{} &  & \multicolumn{1}{c}{} & loss & inf & loss & inf & loss & inf \\
\midrule
\multirow{2}{*}{Bonus.B54.1000} & output & \multirow{2}{*}{3.18} & 3.51 & 28.10 & 4.06 & 36.90 & 4.73 & 40.00 \\
 & program &  & 3.59 & 0.10 & 4.28 & 0.10 & 5.16 & 0.10 \\
\midrule
\multirow{2}{*}{Alstom.Eco.80} & output & \multirow{2}{*}{3.16} & 3.50 & 27.30 & 4.05 & 34.80 & 4.73 & 38.80 \\
 & program &  & 3.55 & 0.30 & 4.20 & 0.30 & 5.01 & 0.30 \\
\midrule
\multirow{2}{*}{Vestas.V27.225} & output & \multirow{2}{*}{3.10} & 3.44 & 11.10 & 3.99 & 24.00 & 4.68 & 29.60 \\
 & program &  & 3.48 & 0.10 & 4.17 & 0.10 & 5.04 & 0.10 \\
\midrule
\multirow{2}{*}{Siemens.SWT.3} & output & \multirow{2}{*}{3.22} & 3.59 & 6.60 & 4.18 & 17.90 & 4.90 & 25.10 \\
 & program &  & 3.59 & 0.00 & 4.23 & 0.00 & 5.02 & 0.00 \\
\midrule
\multirow{2}{*}{Enercon.E92.2300} & output & \multirow{2}{*}{3.20} & 3.56 & 6.00 & 4.13 & 17.30 & 4.84 & 24.60 \\
 & program &  & 3.56 & 0.80 & 4.17 & 0.80 & 4.95 & 0.80 \\
\midrule
\multirow{2}{*}{GE.2.75.103} & output & \multirow{2}{*}{3.08} & 3.43 & 0.00 & 3.99 & 1.80 & 4.68 & 6.00 \\
 & program &  & 3.43 & 0.00 & 4.00 & 0.30 & 4.75 & 0.30 \\
\bottomrule
\end{tabular}
}
{}
\end{table}

\subsection{Private Semidefinite Programming: Maximum Volume Inscribed Ellipsoid}\label{subsec:elliposid}
Let $\mathcal{U}$ be a bounded, non-empty polyhedral set given by a number of linear equalities:
\begin{align*}
    \mathcal{U} = \{x\in\mathbb{R}^{n}\;|\;a_{i}^{\top}x\leqslant b_{i},\;\forall i=1,\dots,m\}.
\end{align*}
The goal is to compute the maximum volume ellipsoid $\mathcal{E}$ contained in $\mathcal{U}$ described by a set
\begin{align*}
    \mathcal{E} = \{x=Yu+z\;|\;\norm{u}\leqslant1\}
\end{align*}
with design parameters $Y\in\mathbb{R}^{n\times n}$ and $z\in \mathbb{R}^{n}$. This problem often arises in the context of robust optimization, where $\mathcal{U}$ is the box uncertainty set and the maximum volume ellipsoid is computed to reduce the conservatism of $\mathcal{U}$ \citep{ben1999robust,cohen2020feature}. According to \citep[Proposition 3.7.1]{ben2011lectures}, the optimal ellipsoid parameters are provided by the solution of the following convex, semidefinite program:
\begin{subequations}\label{MVE:det}
\begin{align}
    \maximize{z,Y}\quad&\text{det}[Y]^{\frac{1}{n}}\label{MVE:det_obj}\\
    \st\quad&Y\succcurlyeq 0,\quad\norm{Ya_{i}}_{2}\leqslant b_{i} - a_{i}^{\top}z, \quad\forall i = 1,\dots, m,
\end{align}
\end{subequations}
where operator $\text{det}[Y]$ denotes the determinant of matrix $Y$. As per \cite[Theorem 3.7.1]{ben2011lectures}, the optimal ellipsoid parameters $z^{\star}$ and $Y^{\star}$ are unique. Hence, changes in set $\mathcal{U}$ are leaking when answering identity queries on the ellipsoid parameters. The differential privacy goal is thus to make adjacent sets $\mathcal{U}$ and $\mathcal{U}'$ statistically similar in perturbed ellipsoid parameters. Using Proposition \ref{prop:OP_IP_dp}, this privacy goal is achieved by either perturbing $z^{\star}$ and $Y^{\star}$ (output perturbation) or by solving program \eqref{MVE:det} on perturbed set $\mathcal{U}$ (input perturbation). In both cases, there is no guarantee that the resulting ellipsoid will be inscribed in the original set. Following the proposed program perturbation strategy from Section \ref{sec:main}, we define linear decision rules for ellipsoid parameters as:
\begin{subequations}\label{Ellipoid_LDR}
\begin{align}
    z(\bs{\zeta}) =\;& \overline{z} + Z\bs{\zeta},\;\overline{z}\in\mathbb{R}^{n}, \;Z\in\mathbb{R}^{n\times3n}\\
    Y(\bs{\zeta}) =\;& \overline{Y} + \tilde{Y}\bs{\zeta},\;\overline{Y}\in\mathbb{R}^{n\times n}, \;\tilde{Y}\in\mathbb{R}^{n\times n\times3n}.
\end{align}
\end{subequations}
Then, the privacy goal is achieved by perturbing the nominal solution $(\overline{z},\overline{Y})$ of the following chance-constrained counterpart of program \eqref{MVE:det}:
\begin{subequations}\label{MVE:sto}
\begin{align}
    \maximize{t,\overline{z},\overline{Y},Z,\tilde{Y}}\quad&t\\
    \st\quad&
    \text{Pr}\left[
    \begin{array}{l}
    t\leqslant\text{det}[\overline{Y} + \tilde{Y}\bs{\zeta}]^{\frac{1}{n}},\\
    \overline{Y} + \tilde{Y}\bs{\zeta}\succcurlyeq 0,\\
    \norm{(\overline{Y} + \tilde{Y}\bs{\zeta})a_{i}}_{2}\leqslant b_{i} - a_{i}^{\top}(\overline{z} + Z\bs{\zeta}), \quad\forall i = 1,\dots, m
    \end{array}
    \right]\geqslant 1-\eta,\label{MVE:sto_cc}\\
    &Z,\tilde{Y}\in \mathcal{X},\label{MVE:sto_identity}
\end{align}
\end{subequations}
where constraint \eqref{MVE:sto_cc} ensures that the ellipsoid with perturbed parameters lies inside the polyhedral set with probability at least $1-\eta\%.$ \textcolor{maincolor}{Unlike previous cases of LP and QP, there is no analytic expression for the expected value of \eqref{MVE:det_obj} under perturbation, and we moved this non-linear term inside the chance constraint using auxiliary variable $t.$ Then, we enforce the entries of chance constraint \eqref{MVE:sto_cc} on the vertices of the rectangular uncertainty set, as explained in Section \ref{subsec:tract_ref}.} Constraint \eqref{MVE:sto_identity} requires the recourse of decision rules \eqref{Ellipoid_LDR} to be independent from the parameters of the polyhedral set. For $n=2$, $\mathcal{X}$ is such that:
$$Z^{\star}=
\begin{bmatrix}
I_{2\times2} & \mathbb{0}_{2\times4}
\end{bmatrix},\quad 
\tilde{Y}^{\star}=
\begin{bmatrix}
\begin{bmatrix}
\mathbb{0}_{2\times2} & I_{2\times2} & \mathbb{0}_{2\times2}
\end{bmatrix} & 
\begin{bmatrix}
\mathbb{0}_{2\times2} & \mathbb{0}_{2\times2} & I_{2\times2}^{\top}
\end{bmatrix}
\end{bmatrix},
$$
and the private identity query is $z(\bs{\zeta}) =\; \overline{z} + [\bs{\zeta}]_{1:2},$ $Y(\bs{\zeta}) =\; \overline{Y} + \begin{bmatrix}
    [\bs{\zeta}]_{3:4} & [\bs{\zeta}]_{5:6}
    \end{bmatrix}.$

\begin{figure}
\FIGURE{
\includegraphics[width=0.95\textwidth]{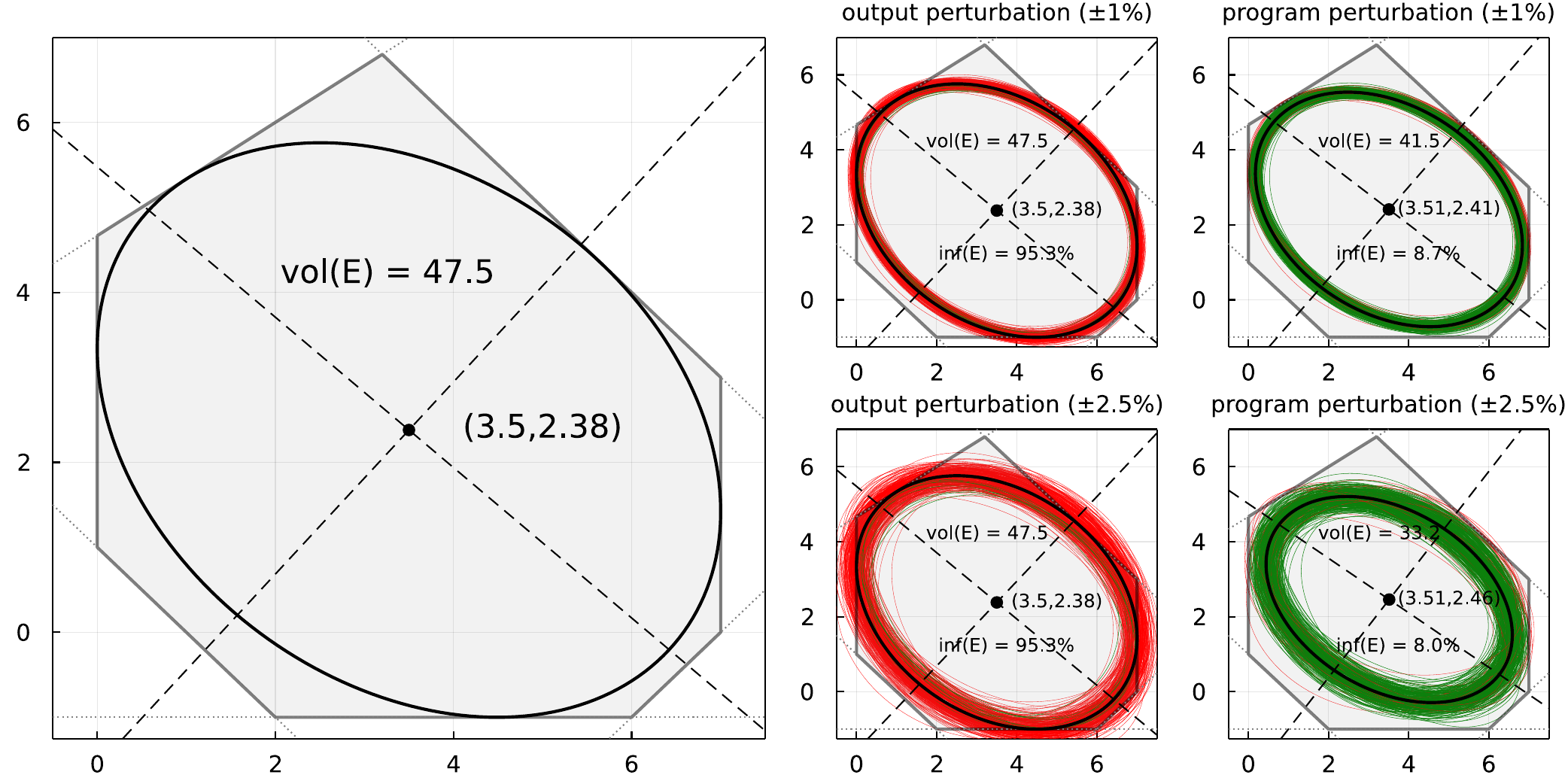}}
{
Privacy-preserving maximum volume inscribed ellipsoid: Deterministic (non-private) solution on the left, and privacy-preserving solution on the right based on output and program perturbation strategies. Perturbing the parameters of the mean ellipsoid (solid black line) results in private ellipsoid models that are outside the polyhedral set (red samples) or inscribed in the set (green samples). The program perturbation produces inscribed ellipsoids with a high probability, yet it is suboptimal w.r.t. the non-private solution in terms of the mean ellipsoid volume.\label{fig:Ellipsoid}} 
{}
\end{figure}

We illustrate the private ellipsoid computation for a polyhedral set depicted in Figure \ref{fig:Ellipsoid} (left). We model a dataset universe where vector $b$ belongs to a range $[b-\gamma,b+\gamma]$ and choose to provide $(1,0.1)-$differential privacy calibrating parameters of perturbation $\bs{\zeta}$ to Gaussian distribution. We require perturbed ellipsoids to remain within the set with a probability of at least $90\%$ $(\eta=10\%)$, \textcolor{maincolor}{and we use the confidence parameter $\beta=10\%$ for Proposition \ref{prop:sample_size_1}, resulting in $S=179$ samples.} 

Figure \ref{fig:Ellipsoid} (right) illustrates the private ellipsoids under the output and program perturbation strategies for $\gamma=1\%$ (top) and $\gamma=2.5\%$ (bottom) of the nominal value of $b$. Observe, the mean ellipsoid under program perturbation features smaller volume than under output perturbation. However, the perturbation of its parameters results in ellipsoids that remain within the polyhedral set with a high probability, revealing the trade-off between the feasibility of the private ellipsoids and their optimality.

\section{Conclusions and perspectives}\label{sec:conclusion}

We proposed a new privacy-preserving stochastic optimization framework for convex optimization with embedded  Laplace and Gaussian differential privacy mechanisms. Specifically, we conditioned optimization outcomes on the realization of privacy-preserving random perturbations, which are accommodated in the feasible region in a least-cost manner using chance-constrained programming. Our framework spans across conic optimization, enabling differentially private linear, quadratic and semidefinite programming in many real-world optimization and machine learning contexts. 

This work has been driven by problems solved in a centralized manner, e.g., when there is a trustworthy centralized agent who collects data and then privately solves an optimization program. Although this setup is standard for many real-life problems, such as energy optimization problems considered in this paper, in the future, data owners may give preferences to decentralized or distributed optimization, so as to prevent data leakage associated with external data transfers and storage. For this future scenario, it is relevant to study how linear decision rules can be optimized in a decentralized/distributed manner, to enable \textcolor{maincolor}{local differential privacy} guarantees while preserving feasibility and optimality of the centralized model. 

\newpage
\bibliographystyle{informs2014} 
\bibliography{references.bib} 

\newpage

\ECSwitch


\ECHead{Supplementary Material}


\section{Proofs of Statements}\label{app:proofs}

\subsection{Proof of Proposition \ref{prop:OP_IP_dp}}
The privacy of the output perturbation is a direct consequence of Theorem \ref{th:Laplace}. In the input perturbation, the perturbed dataset $\tilde{\mathcal{D}}$ is $\varepsilon-$differentially private due to Theorem \ref{th:Laplace}. Since map $x(\cdot)$ is a data-independent transformation of datasets, the input perturbation is also $\varepsilon-$differentially private by post-processing immunity of Theorem \ref{th:post_processing}.

\subsection{Proof of Theorem \ref{th:eps_dp_id}}
We need to show that the probability ratio of observing the same outcome $\hat{x}$ is bounded by $\varepsilon$ as
\begin{subequations}
\begin{align}\label{eq:ratio_1}
    \frac{\text{Pr}[\overline{x}^{\star}(\mathcal{D}) + X(\mathcal{D})^{\star}\bs{\zeta} \in \hat{x}]}{\text{Pr}[\overline{x}^{\star}(\mathcal{D}') + X^{\star}(\mathcal{D}')\bs{\zeta} \in \hat{x}]}\leqslant\text{exp}(\varepsilon).
\end{align}
Due to a query-specific feasible set $\mathcal{X}=\{X|X=I\}$ in program \eqref{problem:SP}, we know that $X^{\star}=I$ in the optimum regardless of the optimization dataset. Hence, the probability ratio $\eqref{eq:ratio_1}$ is equivalent to
\begin{align}
    \frac{\text{Pr}[\overline{x}^{\star}(\mathcal{D}) + \bs{\zeta} \in \hat{x}]}{\text{Pr}[\overline{x}^{\star}(\mathcal{D}') + \bs{\zeta} \in \hat{x}]}\leqslant\text{exp}(\varepsilon).
\end{align}
To show that this inequality holds, similarly to \cite[Theorem 3.6]{dwork2014algorithmic}, we expand the left-hand side as 
\begin{align}
    &\frac{\text{Pr}[\overline{x}^{\star}(\mathcal{D}) + \bs{\zeta} \in \hat{x}]}{\text{Pr}[\overline{x}^{\star}(\mathcal{D}') + \bs{\zeta} \in \hat{x}]} 
    = 
    \frac{\text{Pr}[\bs{\zeta} \in \hat{x} - \overline{x}^{\star}(\mathcal{D}) ]}{\text{Pr}[\bs{\zeta} \in \hat{x} - \overline{x}^{\star}(\mathcal{D}')]} 
    \overset{(\star)}{=}
    \frac{\prod_{i=1}^{n} \text{exp}\left(-\frac{\varepsilon\norm{\hat{x}_{i}-\overline{x}_{i}^{\star}(\mathcal{D})}_{1}}{\Delta_{1}}\right)}{\prod_{i=1}^{n} \text{exp}\left(-\frac{\varepsilon\norm{\hat{x}_{i}-\overline{x}_{i}^{\star}(\mathcal{D}')}_{1}}{\Delta_{1}}\right)}\nonumber\\
    &=\prod_{i=1}^{n}\text{exp}\left(\frac{\varepsilon\norm{\hat{x}_{i} - \overline{x}_{i}^{\star}(\mathcal{D}')}_{1} - \varepsilon\norm{\hat{x}_{i} - \overline{x}_{i}^{\star}(\mathcal{D})}_{1} }{\Delta_{1}}\right)
    \overset{(\pounds)}{\leqslant}
    \prod_{i=1}^{n}\text{exp}\left(\frac{\varepsilon\norm{\overline{x}_{i}^{\star}(\mathcal{D}) - \overline{x}_{i}^{\star}(\mathcal{D}')}_{1}}{\Delta_{1}}\right)\nonumber\\
    &=\text{exp}\left(\frac{\varepsilon\norm{\overline{x}^{\star}(\mathcal{D}) - \overline{x}^{\star}(\mathcal{D}')}_{1}}{\Delta_{1}}\right) 
    \overset{(\S)}{\leqslant} \text{exp}\left(\frac{\varepsilon\Delta_{1}}{\Delta_{1}}\right) = \text{exp}(\varepsilon),
\end{align}
where $(\star)$ is due to definition of the probability density function of the Laplace distribution, $(\pounds)$ follows the reverse inequality of norms, $(\S)$ is from the definition of the worst-case $\ell_{1}-$sensitivity of the identity optimization query to $\alpha$-adjacent datasets.  
\end{subequations}

\subsection{Proof of Theorem \ref{th:eps_dp_sum}}
\begin{subequations}
The steps are similar to the previous proof. We need to show that the probability ratio of observing the same sum query outcome $\hat{x}$ is bounded by $\varepsilon$ as
\begin{align}\label{eq:ratio_2}
    \frac{\text{Pr}[\mathbb{1}^{\top}\left(\overline{x}^{\star}(\mathcal{D}) + X(\mathcal{D})^{\star}\bs{\zeta}\right) \in \hat{x}]}{\text{Pr}[\mathbb{1}^{\top}\left(\overline{x}^{\star}(\mathcal{D}') + X^{\star}(\mathcal{D}')\bs{\zeta}\right) \in \hat{x}]}\leqslant\text{exp}(\varepsilon).
\end{align}
For sum queries, $|\bs{\zeta}|=1$ and we know that program \eqref{problem:SP} includes constraint set $\mathcal{X}=\{X|\mathbb{1}^{\top}X=1\}$ which makes the random query component independent from data, i.e., 
\begin{align}
    \frac{\text{Pr}[\mathbb{1}^{\top}\overline{x}^{\star}(\mathcal{D}) + \bs{\zeta} \in \hat{x}]}{\text{Pr}[\mathbb{1}^{\top}\overline{x}^{\star}(\mathcal{D}') + \bs{\zeta} \in \hat{x}]}\leqslant\text{exp}(\varepsilon).
\end{align}
To prove this inequality, we expand the left-hand side as
\begin{align}
    &\frac{\text{Pr}[\mathbb{1}^{\top}\overline{x}^{\star}(\mathcal{D}) + \bs{\zeta} 
    \in \hat{x}]}{\text{Pr}[\mathbb{1}^{\top}\overline{x}^{\star}(\mathcal{D}') + \bs{\zeta} \in \hat{x}]} 
    = 
    \frac{\text{Pr}[\bs{\zeta} \in \hat{x} - \mathbb{1}^{\top}\overline{x}^{\star}(\mathcal{D})]}{\text{Pr}[\bs{\zeta} \in \hat{x} - \mathbb{1}^{\top}\overline{x}^{\star}(\mathcal{D}')]} 
    \overset{(\star)}{=}
    \frac{\text{exp}\left(\frac{\varepsilon\norm{\hat{x}-\mathbb{1}^{\top}\overline{x}^{\star}(\mathcal{D})}_{1}}{\Delta_{1}}\right)}{\text{exp}\left(\frac{\varepsilon\norm{\hat{x}-\mathbb{1}^{\top}\overline{x}^{\star}(\mathcal{D}')}_{1}}{\Delta_{1}}\right)}\nonumber\\
    &=
    \text{exp}\left(\frac{\varepsilon\norm{\hat{x}-\mathbb{1}^{\top}\overline{x}^{\star}(\mathcal{D}')}_{1} - \varepsilon\norm{\hat{x}-\mathbb{1}^{\top}\overline{x}^{\star}(\mathcal{D})}_{1}}{\Delta_{1}}\right)
    \overset{(\pounds)}{\leqslant}
    \text{exp}\left(\frac{\varepsilon\norm{\mathbb{1}^{\top}\overline{x}^{\star}(\mathcal{D}') - \mathbb{1}^{\top}\overline{x}^{\star}(\mathcal{D}')}_{1}}{\Delta_{1}}\right) \overset{(\S)}{\leqslant} \text{exp}(\varepsilon),\nonumber
\end{align}
where $(\star)$ is due to definition of the probability density function of the Laplace distribution, $(\pounds)$ follows the reverse inequality of norms, $(\S)$ is from the definition of the worst-case $\ell_{1}-$sensitivity of the sum query to $\alpha$-adjacent datasets.  
\end{subequations}

\subsection{Proof of Theorem \ref{th:eps_pdp_id}}

\textcolor{maincolor}{
We analyze the Laplace mechanism with perturbation $\boldsymbol{\zeta}\sim\text{Lap}(t^{\star}/\varepsilon)$. For dataset $\mathcal{D}$, we have
\begin{align*}
    &\text{Pr}\big[\tilde{x}(\mathcal{D})\in\hat{x}\big] = \text{Pr}\big[\overline{x}^{\star}(\mathcal{D})+X^{\star}(\mathcal{D})\boldsymbol{\zeta}\in\hat{x}\big] = \text{Pr}\big[\overline{x}^{\star}(\mathcal{D})+\boldsymbol{\zeta}\in\hat{x}\big] \\
    &= \text{Pr}\big[\boldsymbol{\zeta}\in\hat{x}-\overline{x}^{\star}(\mathcal{D})\big] \propto
    \text{exp}\left(-\frac{\varepsilon\norm{\hat{x}-\overline{x}^{\star}(\mathcal{D})}}{t^{\star}}\right)
\end{align*}
Likewise, for the adjacent dataset $\mathcal{D'}$ we have 
\begin{align*}
    &\text{Pr}\big[\tilde{x}(\mathcal{D}')\in\hat{x}\big] \propto
    \text{exp}\left(-\frac{\varepsilon\norm{\hat{x}-\overline{x}^{\star}(\mathcal{D}')}}{t^{\star}}\right)
\end{align*}
The log-ratio of probabilities is then
\begin{align*}
    \left|\text{ln}\left(\frac{\text{Pr}\big[\tilde{x}(\mathcal{D}\textcolor{white}{'})\in\hat{x}\big]}{\text{Pr}\big[\tilde{x}(\mathcal{D}')\in\hat{x}\big]}\right)\right|
    &=
    \left|\text{ln}\left(\text{exp}\left(-\frac{\varepsilon\norm{\hat{x}-\overline{x}^{\star}(\mathcal{D})}}{t^{\star}}\right)\right) - \text{ln}\left(\text{exp}\left(-\frac{\varepsilon\norm{\hat{x}-\overline{x}^{\star}(\mathcal{D}')}}{t^{\star}}\right)\right)\right|\\
    &=\left|
    \frac{\varepsilon}{t^{\star}}\norm{\hat{x}-\overline{x}^{\star}(\mathcal{D}')} - 
    \frac{\varepsilon}{t^{\star}}\norm{\hat{x}-\overline{x}^{\star}(\mathcal{D})}
    \right|\\
    &\leqslant\frac{\varepsilon}{t_{1}^{\star}}\norm{x(\mathcal{D}) - x(\mathcal{D}')}
\end{align*}
Hence, when $\norm{x(\mathcal{D}) - x(\mathcal{D}')}\leqslant t^{\star}$,  $\tilde{x}(\mathcal{D})$ satisfies $\varepsilon-$differential privacy. Since $t^{\star}$ solves the chance-constrained program \eqref{prog:sens_cc_program}, we know that for a random draw $\mathcal{D},\mathcal{D}'\sim\mathbb{D}$, the inequality $\norm{x(\mathcal{D}) - x(\mathcal{D}')}\leqslant t^{\star}$ holds with probability at least $1-\gamma$, i.e., $\tilde{x}(\mathcal{D})$ satisfies $\varepsilon-$differential privacy with at least $1-\gamma$ probability. Therefore, $\tilde{x}(\mathcal{D})$ satisfies $(\varepsilon,\gamma)-$probabilistic differential privacy.}

\subsection{Proof of Theorem \ref{th:eps_delta_dp_id}}
Similarly to the proofs of Theorems \ref{th:eps_dp_id}, the random components of the identity query must be shown to be independent from datasets $\mathcal{D}$ and $\mathcal{D}'$ using the query specific feasibility set $\mathcal{X}$. The remainder follows the same steps as in the proof in \cite[Apendix A]{dwork2014algorithmic}, where the function of interest $f(d)$ is substituted with $\tilde{x}(\mathcal{D})$ for the identity query and with $\mathbb{1}^{\top}\tilde{x}(\mathcal{D})$ for the sum query, and function sensitivity $\Delta f$ must be equal to $\Delta_{2}$. 

\textcolor{maincolor}{
\section{Computing the worst-case solution sensitivity to adjacent datasets}
While we use optimization problem \eqref{prog:sens_sample} to streamline the minimum sample size requirement, we don't have to solve \eqref{prog:sens_sample} as an optimization problem. Here, we put forth a simple Algorithm \ref{alg:sens} to compute the sensitivity. In this algorithm, Step 1 repeats $S$ times to statistically bound the sensitivity. The while loop at Step 3 repeats until a pair of $\alpha-$adjacent datasets is drawn from dataset universe $\mathbb{D}$ at Step 4 and the sensitivity to this dataset pair is evaluated at Step 6. The last step 10 takes the maximal sensitivity observed over iterations. }
\begin{algorithm}[h!]
\SetKwInOut{Input}{Input}
\SetKwInOut{Output}{Output}
\caption{Sensitivity computation.\label{alg:singlepm}}\label{alg:sens}
\Input{%
		$\alpha$ -- adjacency, $p$ -- order, $S$ -- sample size
	  }
\Output{%
		\xvbox{4mm}{$\Delta_{p}$} -- sensitivity approximation
	   }
  \BlankLine 
  \For{$s=1,\dots,S$ }{
  initialize $\mathcal{D}$ and $\mathcal{D}'$\\
    \While{$\norm{\mathcal{D}-\mathcal{D}'}>\alpha$}{
        sample $\mathcal{D},\mathcal{D}'\sim\mathbb{D}$\\
            \If{ $\norm{\mathcal{D}-\mathcal{D}'}\leqslant\alpha$}{
            $\delta_{s}=\norm{x(\mathcal{D}) - x(\mathcal{D}')}_{p}$ 
            }
    }
  }
  \xvbox{4mm}{$\Delta_{p}$} $\leftarrow$ $\text{max}\{\delta_{1},\dots,\delta_{S}\}$ 
\end{algorithm}

\section{Conic Representation of Problems of Interest}\label{app:conic_representations}

In this section, we demonstrate how some problems of interest are indeed representable by the reference conic program \eqref{problem:base}, so that the results of Section \ref{sec:main} apply to applications in Section \ref{sec:application}.  

\begin{problem}[Optimal power flow]\normalfont A linear optimization program 
\begin{subequations}
\begin{align}
    \minimize{x}\quad&c^{\top}x\\
    \st\quad&\mathbb{1}^{\top}(x-d)=0,\\
    &|F(x-d)|\leqslant f^{\text{max}},\\
    &x^{\text{min}}\leqslant x\leqslant x^{\text{max}},
\end{align}
\end{subequations}
with variable $x\in\mathbb{R}^{N}$ and parameters $c,d,x^{\text{min}},x^{\text{max}}\in\mathbb{R}^{N},F\in\mathbb{R}^{E\times N}$ and $f^{\text{max}}\in\mathbb{R}^{E}$ is presentable by conic program \eqref{problem:base} when cone 
$\mathcal{K}=\mathbb{R}_{+}^{2+2N+2E}$, and the dataset is
\begin{align*}
\begingroup
\renewcommand*{\arraystretch}{1.15}
\mathcal{D} = 
\left(
\begin{bmatrix}
-\mathbb{1}^{\top}\\
\textcolor{white}{-}\mathbb{1}^{\top}\\
\textcolor{black}{-}F\\
\textcolor{white}{-}F\\
\textcolor{black}{-}I\\
\textcolor{white}{-}I\\
\end{bmatrix}, 
\begin{bmatrix}
\textcolor{white}{-}\mathbb{1}^{\top}d\\
\textcolor{black}{-}\mathbb{1}^{\top}d\\
\textcolor{white}{-}f^{\text{max}} - Fd\\
\textcolor{white}{-}f^{\text{max}} - Fd\\
\textcolor{white}{-}x^{\text{max}}\\
\textcolor{black}{-}x^{\text{min}}\\
\end{bmatrix},
c
\right).
\endgroup
\end{align*}
\end{problem}

\begin{problem}[Support vector machines]\normalfont
A quadratic optimization program
\begin{subequations}
\begin{align}
    \minimize{w,b,z}\quad&\lambda\norm{w}^{2} + \textstyle\frac{1}{m}\mathbb{1}^{\top}z\\
    \st\quad&y_{i}(w^{\top}x_{i}-b)\geqslant1-z_{i},\quad z_{i}\geqslant0,\quad\forall i=1,\dots,m,
\end{align}
\end{subequations}
with variables $w\in\mathbb{R}^{n},b\in\mathbb{R},z\in\mathbb{R}^{m}$ and parameters $y\in\mathbb{R}^{m}$ and $x_{1},\dots,x_{m}\in\mathbb{R}^{n}$ can be equivalently rewritten using an auxiliary scalar variable $t$ as
\begin{subequations}
\begin{align}
    \minimize{t,w,b,z}\quad&\lambda t + \textstyle\frac{1}{m}\mathbb{1}^{\top}z\\
    \st\quad& \norm{w}^{2}\leqslant t,\\
    &\text{diag}[y] 
    \left(\begin{bmatrix}
    x_{1}\\\vdots\\x_{m}
    \end{bmatrix}
    w - b\right) + z \geqslant\mathbb{1},\\
    &z\geqslant\mathbb{0},
\end{align}
\end{subequations}
which is presentable by conic program \eqref{problem:base} when we have
\begin{align*}
x = 
\begin{bmatrix}
t\\
w\\
b\\
z
\end{bmatrix},\;
\mathcal{D}=
\left(
\begin{bmatrix}
-1 & \vertbar & \vertbar & \vertbar\\
\vertbar & -I & \vertbar & \vertbar\\
\cmidrule(lr){1-4}
\vertbar & -\text{diag}[y] 
\begin{bmatrix}
x_{1}\\\vdots\\x_{m}
\end{bmatrix}
& \vertbar & \vertbar\\
\vertbar & \vertbar & - \text{diag}[y] & \vertbar \\
\vertbar & \vertbar & \vertbar & \textcolor{white}{-}I \\
\cmidrule(lr){1-4}
\vertbar & \vertbar & \vertbar & \textcolor{black}{-}I 
\end{bmatrix},
\begin{bmatrix}
\vertbar\\
\cmidrule(lr){1-1}
\textcolor{black}{-}\mathbb{1}\\
\cmidrule(lr){1-1}
\textcolor{white}{-}\mathbb{0}
\end{bmatrix},\;
\begin{bmatrix}
1\\
\textstyle\frac{1}{m}\cdot\mathbb{1}\\
\vertbar\\
\vertbar
\end{bmatrix}
\right),\;
\mathcal{K}=\mathcal{Q}_{n+1}^{r}\times\mathbb{R}_{+}^{m}\times\mathbb{R}_{+}^{m},
\end{align*}
where $\mathcal{Q}_{n+1}^{r}$ is a rotated second-order cone, blank spaces are zeros and horizontal lines are cone boundaries. 
\end{problem}

\begin{problem}[Monotonic regression]\normalfont A quadratic optimization program 
\begin{subequations}
\begin{align}
\minimize{w}\quad&\sum_{i=1}^{n}\norm{y_{i} - w^{\top}\varphi(x_{i})}^2 + \lambda\norm{w}^2\\
\st\quad&Cw\geqslant\mathbb{0},
\end{align}
\end{subequations}
with variable $w\in\mathbb{R}^{m}$, parameters $y\in\mathbb{R}^{n},x_{1},\dots,x_{n}\in\mathbb{R}^{k},\lambda\in\mathbb{R}$ and $C\in\mathbb{R}^{p\times m}$, and some basis function $\varphi:\mathbb{R}^{k}\mapsto\mathbb{R}^{m}$, takes the following equivalent conic form
\begin{subequations}
\begin{align}
\minimize{u,\nu,w}\quad&u + \lambda\nu\\
\st\quad& \norm{y - 
\begin{bmatrix}
\varphi(x_{1}) \\
\vdots \\
\varphi(x_{n}) 
\end{bmatrix}
w}^{2}\leqslant u,\\
& \norm{w}^{2}\leqslant \nu,\\
&Cw\geqslant\mathbb{0},
\end{align}
\end{subequations}
using auxiliary variables $u,\nu\in\mathbb{R}$. This problem is representable by conic program \eqref{problem:base} when
\begin{align*}
x = 
\begin{bmatrix}
u\\
\nu\\
w
\end{bmatrix},\;
\mathcal{D}=\left(
\begin{bmatrix}
1 & \vertbar & \vertbar \\
\vertbar & \vertbar & 
\begin{bmatrix}
\varphi(x_{1}) \\
\vdots \\
\varphi(x_{n}) 
\end{bmatrix} \\
\cmidrule(lr){1-3}
\vertbar & 1 & \vertbar \\
\vertbar & \vertbar & -I \\
\cmidrule(lr){1-3}
\vertbar & \vertbar & -C \\
\end{bmatrix},
\begin{bmatrix}
0\\
y\\
\cmidrule(lr){1-1}
\vertbar\\
\cmidrule(lr){1-1}
\vertbar
\end{bmatrix},
\begin{bmatrix}
1\\
\lambda\\
\vertbar
\end{bmatrix}
\right),\; 
\mathcal{K}=\mathcal{Q}_{n+1}^{r}\times \mathcal{Q}_{m+1}^{r}\times \mathbb{R}_{+}^{p},
\end{align*}
where horizontal lines draw cones boundaries among two rotated second-order cones and one non-negative orthant. 
\end{problem}

\begin{problem}[Maximum-volume inscribed ellipsoid]\normalfont
Consider a semidefinite program 
\begin{subequations}\label{sdp_example}
\begin{align}
    \maximize{z,Y}\quad&\text{det}[Y]^{\frac{1}{n}}\label{obj_det}\\
    \st\quad&\norm{Ya_{i}}_{2}\leqslant b_{i} - a_{i}^{\top}z, \quad\forall i = 1,\dots, m,\quad Y\succcurlyeq 0
\end{align}
\end{subequations}
with variables $z\in\mathbb{R}^{n},Y\in\mathbb{R}^{n\times n}$ and parameters $b\in\mathbb{R}^{m},a_{1},\dots,a_{m}\in\mathbb{R}^{n}$. According to \citet[Section 3.2]{ben2011lectures}, for a fixed $Y^{\star}$, evaluating function \eqref{obj_det} amounts to solving a conic program
\begin{subequations}\label{det_ref}
\begin{align}
    \maximize{t,K}\quad&t\\
    \st\quad
    &\begin{bmatrix}
    t\\
    (I\circ K)\mathbb{1}
    \end{bmatrix}
    \in\mathcal{G}_{n+1}^{\frac{1}{n}},\quad 
    \begin{bmatrix}
    Y^{\star} & K\\
    K^{\top} & I\circ K
    \end{bmatrix}
    \in\mathcal{S}^{+}_{2n},
\end{align}
\end{subequations}
with auxiliary variables $t\in\mathbb{R}$ and $K\in\mathbb{R}^{n\times n}$, one power cone $\mathcal{G}_{n+1}^{\frac{1}{n}}$ and one positive semidefinite cone $\mathcal{S}^{+}_{2n}$. Substituting \eqref{obj_det} with \eqref{det_ref} leads to a conic reformulation of semidefinite program \eqref{sdp_example}. This problem is brought to the standard form in \eqref{problem:base} when 
\begin{align*}
    x = \begin{bmatrix*}[r]
    t\\
    z\\
    \text{vec}[Y]\\
    \text{vec}[K^{\top}]\\
    \text{vec}[K^{\textcolor{white}{\top}}]\\
    \text{vec}[I\circ K^{\textcolor{white}{\top}}]
    \end{bmatrix*},\quad 
    \mathcal{K}= \mathcal{S}^{+}_{n} \times \mathcal{Q}_{n+1}^{1} \times \dots \times \mathcal{Q}_{n+1}^{m} \times \mathcal{G}_{n+1}^{\frac{1}{n}} \times \mathcal{S}^{+}_{2k},
\end{align*}
and the dataset is
\begin{align*}
\mathcal{D} = \left(
\begin{bmatrix}
\vertbar &\vertbar & -I_{n} & \vertbar  &\quad \vertbar & \vertbar \\
\cmidrule(lr){1-6}
\vertbar& a_{1} & \vertbar & \vertbar & \quad \vertbar& \vertbar \\
\vertbar& \vertbar& -a_{1}\otimes I_{n} & \vertbar & \quad \vertbar & \vertbar \\
\cmidrule(lr){1-6} 
& & & \dots & &  \\
\cmidrule(lr){1-6}
\vertbar& a_{m} & \vertbar & \vertbar & \quad \vertbar& \vertbar \\
\vertbar& \vertbar& -a_{m}\otimes I_{n} & \vertbar & \quad \vertbar & \vertbar \\
\cmidrule(lr){1-6} 
1 & \vertbar& \vertbar & \vertbar & \quad \vertbar & \vertbar \\
\vertbar& \vertbar& \vertbar & \vertbar &\quad \vertbar & [e_{1}e_{1}^{\top}\dots e_{n}e_{n}^{\top}]\\
\cmidrule(lr){1-6}
\vertbar& \vertbar& \multicolumn{4}{c}{I_{n}\;\otimes\;
\begin{bmatrix}
I_{n}\;\otimes\; [I_{n}\;\;\mathbb{0}_{n}]\\
I_{n}\;\otimes\; [\mathbb{0}_{n}\;\;I_{n}]
\end{bmatrix}
}
\end{bmatrix}, 
\begin{bmatrix} 
\\
\cmidrule(lr){1-1}
b_{1}\\
\\
\cmidrule(lr){1-1}
\dots \\
\cmidrule(lr){1-1}
b_{m}\\
\\
\cmidrule(lr){1-1}
\\
\\
\cmidrule(lr){1-1}\\
\begin{matrix}
\\
\end{matrix}
\end{bmatrix},
\begin{bmatrix}
1\\
\vertbar\\
\vertbar\\
\vertbar\\
\vertbar\\
\vertbar
\end{bmatrix}
\right),
\end{align*}
where $\text{vec}[\cdot]$ is a matrix vectorization operator, $\otimes$ is the Kronecker product, $e_{i}$ is the basis vector with 1 at position $i$ and the rest being zeros, and horizontal lines are cone separation boundaries. 
\end{problem}

%
%
%





\end{document}